\documentclass[a4paper, reqno, 12pt]{amsart}

\usepackage{geometry}       

\usepackage{float}
\usepackage{bm}
\usepackage{fullpage,xcolor}
\usepackage[mathscr]{euscript}
\usepackage[all]{xy}
\usepackage{epsfig}
\usepackage{amsfonts}
\usepackage{mathptmx}

\usepackage{graphicx}
\usepackage{amssymb} 
\usepackage{amsmath}
\usepackage{amsthm}
\usepackage{mathrsfs}
\usepackage{epstopdf}
\usepackage{url}
\usepackage[msc-links,alphabetic]{amsrefs}
\usepackage{tikz}

\usepackage{color}
\makeatletter

\@addtoreset{equation}{section}
\makeatother 

\theoremstyle{plain}
\newtheorem{theorem}{Theorem}[section]

\newtheorem{proposition}[theorem]{Proposition}
\newtheorem{lemma}[theorem]{Lemma}

\theoremstyle{definition}
\newtheorem{definition}[theorem]{Definition}

\newtheorem{remark}[theorem]{Remark}

\usepackage{latexsym}

\title[Equivalence of Doubly Periodic Tangles]
 {Equivalence of Doubly Periodic Tangles}

\author{Ioannis Diamantis}
\address{Department of Data Analytics and Digitalisation,
Maastricht University, School of Business and Economics,
P.O.Box 616, 6200 MD, Maastricht,
The Netherlands.}
\email{i.diamantis@maastrichtuniversity.nl}

\author{Sofia Lambropoulou}
\address{School of Applied Mathematical and Physical Sciences, National Technical University of Athens, Zografou campus, GR-15780 Athens, Greece.}
\email{sofia@math.ntua.gr}
\urladdr{http://www.math.ntua.gr/~sofia}

\author{Sonia Mahmoudi}
\address{Advanced Institute for Materials Research (WPI-AIMR), Tohoku University, 2-1-1 Katahira, Aoba-ku, Sendai, Miyagi 980-8577, Japan.} 
\email{sonia.mahmoudi@tohoku.ac.jp}
\address{RIKEN, Center for Interdisciplinary Theoretical and Mathematical Sciences (iTHEMS), 2-1, Hirosawa, Wako, Saitama 351-0198, Japan.}



\date{}

\thanks{This work was supported by JSPS KAKENHI Grant-in-Aid for Early-Career Scientists (Grant Number 25K17246) and the Daiichi-Sankyo ``Habataku'' Support Program for the Next Generation of Researchers 2025. The authors also acknowledge financial support from the Tohoku Forum for Creativity of Tohoku University, the Japan Tourism Agency (MICE), and the Sendai Tourism, Convention and International Association (SenTIA) for supporting the program \textit{The Theory of Periodic Tangles and Their Interdisciplinary Applications}, during which a significant part of this research was developed. This work was also supported by the National Technical University of Athens and the Maastricht University.}

\begin{document}

\setcounter{section}{-1}

\begin{abstract} 
Doubly periodic tangles, or \textit{DP tangles}, are embeddings of curves in the thickened plane that are periodically repeated in two directions. They are  defined as universal covers of their generating cells, the {\it flat motifs}, which represent knots and links in the thickened torus, and which can be chosen in infinitely many ways. DP tangles are used in modeling materials and physical systems of entangled filaments. In this paper we establish the complete mathematical framework of the topological theory of DP tangles. We present an exhaustive analysis of DP tangle isotopies. These are distinguished in local isotopies and global isotopies. Our analysis  yields the characterization of DP isotopy as an equivalence relation on the level of their (flat) motifs, called \textit{DP tangle equivalence}. Along the way we also discuss motif minimality. We further generalize our results to other diagrammatic categories, namely  framed, virtual, welded, singular,  pseudo, tied and  bonded DP tangles, which could be used in novel applications.
\end{abstract}

\maketitle

\section{Introduction}\label{sec:0}

Doubly periodic tangles, \textit{DP tangles}, are complex entanglements of curves embedded in the thickened plane $\mathbb{E}^2 \times I$ that are periodically repeated in two transversal directions. So, a DP tangle can be defined as the lift to the universal cover $\mathbb{E}^2 \times I$ of a knot or link   in the thickened torus, $T^2 \times I$, called \textit{motif}. 

Periodic tangles are appropriate for modeling and studying materials and physical systems of entangled filaments in various scales, such as polymer melts \cites{Eleni1, Eleni2, Eleni3}, textiles \cites{Sonia1,Sonia2,Isonemal1,Isonemal2,Isonemal3, SatinsTwills,Sabetta}, cosmic filaments \cites{Bond,Hong1,Hong2}, among others. A better understanding of their geometry and topology, often associated to some physical and mechanical properties, could allow the prediction of some of their functions. Following this motivation, Evans, Hyde et al. describe and enumerate periodic tangles using graphs and tilings of the Euclidean and hyperbolic planes \cites{Evans1,Evans2,Evans3,Hyde}. M. O'Keeffe et al. study periodic tangles based on symmetries assumptions, considering PL embeddings with sticks to model structures in molecular chemistry \cites{Yaghi,Treacy}. The first effort toward establishing topological equivalence of DP tangles by means of their motifs was made by Grishanov et al. in \cite{Grishanov1}, mainly in the context of classification of textiles. Yet, there is no universal mathematical study of DP tangles and there are many open questions. 

In this paper we establish the complete mathematical framework of the topological theory of DP tangles. We call two DP tangles \textit{DP isotopic} if they are related via an ambient isotopy of the thickened plane that preserves double periodicity. DP tangle isotopy can be reduced to the diagrammatic level where isotopies are discretized as sequences of periodic moves on \textit{DP diagrams}, generalizing the classical Reidemeister moves. DP tangle isotopy differs from the standard isotopy of classical knots and links in the thickened torus in that it may require an infinite sequence of $\Delta$-moves, which are local isotopy moves. By definition, a DP tangle comes equipped with a non-unique periodic integer lattice, which is determined by a choice of a covering map. 

The aim of the paper is to translate DP tangle isotopy into an equivalence at the level of (flat) motif diagrams, independently of the choice of lattice and basis, unlike previous studies. We will refer to this equivalence relation as \textit{DP tangle equivalence}. For this, we analyze all DP tangle isotopies, which we distinguish in  local and global isotopies. 

Local DP isotopies are analyzed in Section~\ref{sec:local_isotopy}, by exploiting the theory of mixed links \cite{LR1} in order to derive a complete set of local isotopy moves between (flat) motif diagrams. A mixed link represents a link in the thickened torus and comprises the Hopf link as a fixed sublink representing the thickened torus and a `moving' sublink representing the link in the thickened torus. The mixed link point of view is also adopted in \cite[Section 2.1]{Morton}, where a mixed link diagram representing a motif of a DP tangle is called a \textit{kernel}. In \cite{Morton} a motif of a textile structure is modeled by a mixed link as a tool for detecting the existence of null-homotopic components, without addressing isotopy. 

Global DP isotopies arise from affine transformations of the plane. These 
are addressed in Sections~\ref{sec:shiftequivalence},~\ref{sec:scaleequivalence} and~\ref{sec:dehnequivalence}, treating separately translations, scalings and shearings. In each section our strategy is to analyze the following cases: transformation of the lattice keeping the DP diagram point-wise fixed, deformation of the DP diagram keeping the lattice point-wise fixed  and  simultaneous transformation of the DP diagram and its lattice. We prove that the first and second cases are reverse operations, by arguing that an affine transformation applied on a DP diagram can be realized by an infinitum of sequential surface isotopies. We further show that a simultaneous transformation is visible on the DP diagram level but in most cases is not visible on the level of (flat) motifs. We also distinguish between integral and non-integral transformations, highlighting the cases where the DP diagram and/or the lattice remain set-wise fixed, such as in the case of an integral translation. 

More precisely, in Section~\ref{sec:shiftequivalence} we show that a \textit{translation} corresponds to a \textit{shift} of the longitude-meridian pair of the torus or, equivalently, by a specific finite sequence of surface isotopies (see Proposition~\ref{prop:shift-equivalent}).  In Section~\ref{sec:scaleequivalence} we show that  \textit{scaling} is related to \textit{torus inflations or deflations} and/or to \textit{scale equivalence}, which arises from different finite covering maps (see Proposition~\ref{prop:scale-equivalent}). In Section~\ref{sec:dehnequivalence} we show that \textit{shearing} is related to \textit{Dehn twist} of the torus or to an \textit{admissible non-$\mathbb{Z}$ shearing} of the flat motif (see Proposition~\ref{prop:dehn-equivalent}). Moreover, in Section~\ref{sec:unimodular} we investigate further the class of unimodular linear transformations, which are all generated by shearings, according to their classification and their effect on DP diagrams, covering thus all area-preserving linear transformations, including rotations.  The above culminate to the following main result  of the paper (Theorem~\ref{th:equivalence} in Section~\ref{sec:DPequivalent}):

\smallbreak

\noindent \textbf{Theorem (DP tangle equivalence)}.
\textit{Two DP tangles are DP isotopic if and only if two (flat) motif diagrams of theirs are related by a finite sequence of surface isotopy moves, Reidemeister moves, shift equivalence, scale equivalence, torus inflation/deflation, and shear equivalence.}

\smallbreak

The definition of equivalence in \cite{Grishanov1} takes into consideration the different transformations mentioned above, while considering only \textit{minimal (flat) motif diagrams}, namely   motif diagrams created by the quotient of the DP diagram by a \textit{fixed maximal point lattice}. Our result above extends to the  set of all possible  motifs of a DP diagram. 

Translating DP tangle isotopy to the level of (flat) motifs is useful for constructing topological invariants for DP tangles, which are mathematical tools for their classification. See for example \cites{Kurlin, DLM, Sonia1, Sonia2, Grishanov1, Grishanov.part1, Grishanov.part2, Grishanov.Vassiliev1, Grishanov.Vassiliev2,Sonia3, Morton, Eleni1, Eleni2, Eleni3}. These invariants give the same values on all minimal motifs of a given DP tangle, but most of them depend on the minimality assumption. Yet, finding a minimal motif is known to be a very complex problem, which for now has only been solved only for a particular class DP weaves \cite{Sonia1}. In Subsection~\ref{sec:minimal-motif} we introduce the notion of  minimal motif, while in Subsection~\ref{sec:minimalmotif} we highlight the subtlety of finding a minimal motif.

In the last Section~\ref{sec:othersettings}, we generalize our results to other diagrammatic settings. We discuss DP tangle equivalence in the settings of \textit{regular} and \textit{framed} isotopies, as well as \textit{virtual}, \textit{welded}, \textit{singular}, \textit{pseudo}, \textit{tied} and \textit{bonded} DP tangle equivalence. See Theorems~\ref{th:regularframedequivalence}, ~\ref{th:virtualweldedequivalence}, ~\ref{th:singularequivalence}, 

\noindent \ref{th:pseudoequivalence}, and ~\ref{th:tiedbondedequivalence}, analogues of Theorem~\ref{th:equivalence}. The detailed analysis for DP tangles applies equally to any diagrammatic category. Thus, for the study of DP tangles related to any of the above topological settings, we only need to adapt our analysis in the context of motif isotopy. It is worth pointing out a subtlety in the scale equivalence of pseudo DP tangles, due to the fact that for a pre-crossing we do not know whether it an over or an under crossing, cf. \cite{DLM-pseudo}.

The topological models of DP tangles in the different settings  can find applications in areas of science such as, molecular chemistry, textiles
and meta-materials. For example, virtual DP tangles are potentially interesting in materials science, where the prevention of friction between strands is desirable. Pseudo DP tangles \cite{DLM-pseudo} can model worn textiles, where distinguishing the relative positions of two entangled strands may not be obvious. Singular DP equivalence gives the possibility to extend finite type invariants of knots and links to DP tangles. Finally,  bonded DP tangles can be relevant in molecular chemistry.

Concluding, Theorem~\ref{th:equivalence}, along with its analogues related to different diagrammatic settings,  will serve as a foundation for future works on the topological study of DP tangles, as for example in the development of new invariants, that advance their classification problem, including our on-going work \cite{MDL}, as well as for potential novel applications.

\section{Topological set-up for DP tangles and motifs}\label{sec:setup}

In this section we introduce the notions of DP tangles and their generating motifs. 

\subsection{DP tangles, motifs and diagrams} \label{sec:DPtangle}

Let $\tau$ be a knot or link in the thickened torus $T^2 \times I$, where $I=[0,1]$ the unit interval. A \textit{knot} is an embedding of a circle into the thickened torus, while a \textit{link} is an embedding of a finite collection of disjoint circles. We shall say `links' for both knots and links. We equip  $T^2$  with a longitude-meridian pair $(l,m)$. Let further $\mathbb{E}^2$ denote the Euclidean plane and let $B=(u,v)$ be a basis of $\mathbb{E}^2$. We consider the covering map 
$$
\rho: \mathbb{E}^2 \rightarrow{} T^2,
$$
\noindent that assigns the longitude $l$ of $T^2$ to $u$ and  the meridian $m$ of $T^2$ to $v$. Note that the covering map $\rho$ extends trivially to a covering map (also denoted) $\tilde{\rho}$: $ \mathbb{E}^2 \times I \rightarrow{} T^2 \times I$.

\begin{definition}\label{def:DP tangle}
Let $\tau$ be a link in $T^2 \times I$, with $T^2$ equipped with the a longitude-meridian pair $(l,m)$. 
The lift $\tau_{\infty}$ of the pair $\big(\tau, (l,m)\big)$ under the covering map $\tilde{\rho}$ is called a \textit{doubly periodic tangle}, or \textit{DP tangle}. The pair $\big(\tau,(l,m)\big)$, or simply $\tau$, is said to be a \textit{motif} of $\tau_{\infty}$. The diagram of the motif $\tau$ in $T^2 \times \{0\}$, denoted by $\big(d,(l,m)\big)$, or simply $d$, is called a \textit{motif diagram}, and the lift of $d$ under $\rho$, denoted by $d_{\infty}$, is called a \textit{doubly periodic diagram}, or \textit{DP diagram}. 
\end{definition}

\begin{figure}[H]
\centerline{\includegraphics[width=5in]{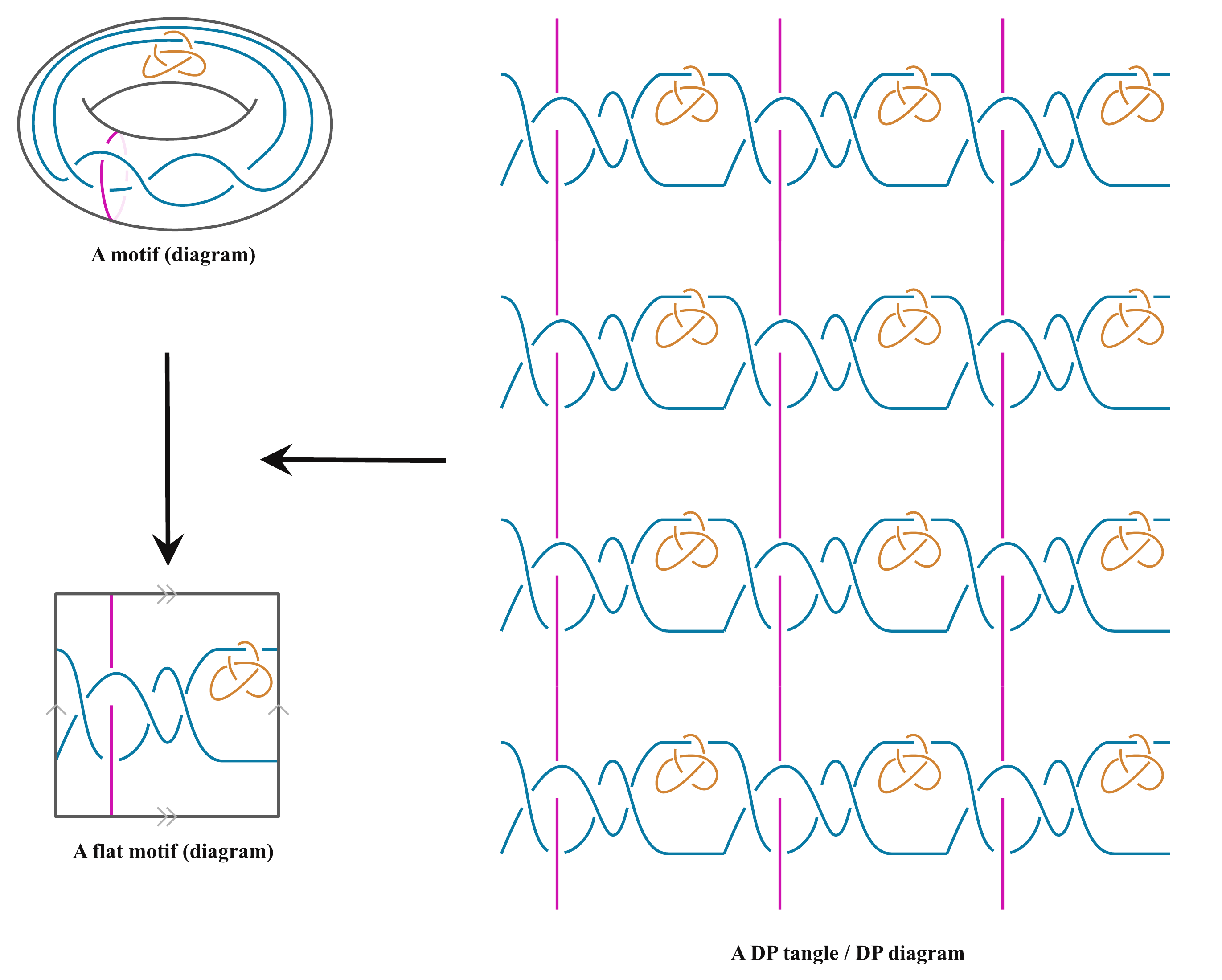}}
\caption{\label{DP-tangle} A motif (diagram), a corresponding flat motif (diagram) and its DP tangle (DP diagram).}
\end{figure}

An example of the above notions is illustrated in Figure~\ref{DP-tangle}. So, a DP tangle comes equipped with a basis of $\mathbb{E}^2$ and a choice of a longitude-meridian pair $(l,m)$ for $T^2$. By a general position argument, $m$ and $l$ do not intersect crossings of $d$ and no arc of $d$ intersects the double point formed by the intersection of $l$ and $m$. Note that in \cite{Grishanov1} a DP tangle is referred to as \textit{doubly periodic structure}, in \cites{Grishanov.part1, Grishanov.part2} as a \textit{2-structure}, while in \cite{Morton} as a \textit{fabric}.

\subsection{Lattices and flat motifs} \label{sec:lattice}

Let $d_\infty$  be a fixed DP diagram. By definition, $d_\infty$ is invariant under translational symmetry along the two linearly independent vectors $u$ and $v$ of the basis $B = (u,v)$ of $\mathbb{E}^2$. The set of points 
\[\Lambda (u,v) = \{xu + yv \ | \ x,y \in \mathbb{Z}\}
\]
generated by the basis $B$ of $\mathbb{E}^2$, defines a \textit{periodic lattice} for the DP diagram $d_\infty$, and by abuse of language of the corresponding DP tangle $\tau_{\infty}$, by trivial extension along the $I$-direction. This implies that the DP tangle $\tau_{\infty}$ arising from the motif $\big(\tau, (l,m)\big)$ comes naturally equipped with the lattice $\Lambda(u,v)$. So we shall  use the notation $(\tau_{\infty},\Lambda)$ and $(d_{\infty},\Lambda)$. On the other hand, starting from $(\tau_{\infty},\Lambda)$ and $(d_{\infty},\Lambda)$, a motif is defined as  $\big(\tau, (l,m)\big) = \tau_{\infty} / \Lambda$ and a motif diagram as $\big(d, (l,m)\big) = d_{\infty} / \Lambda$.

\smallbreak
However, there are infinitely many other pairs  of linearly independent vectors of the plane  giving rise to translational symmetry for the fixed DP diagram  $d_\infty$. Let $\mathcal{T}=\{t \in \mathbb{R}^2\, |\, d_\infty+t = d_\infty\}$ be the translation group of all the translation vectors keeping $d_\infty$ invariant set-wise. It is well-known that $\mathcal{T}$ is isomorphic to $\mathbb{Z}^2$, so there exists a pair of primitive vectors  in $\mathcal{T}$, say $u_0$ and $v_0$, generating $\mathcal{T}$. The set of points
$$
\Lambda_0 (u_0,v_0) = \{xu_0 + yv_0\, |\, x,y \in \mathbb{Z}\} \cong \mathcal{T} \cong \mathbb{Z}^2
$$ 
generated by the two primitive vectors $u_0$ and $v_0$, defines the \textit{maximal periodic lattice} for the DP diagram $d_\infty$ and for the corresponding DP tangle $\tau_{\infty}$. Therefore, any pair of translation vectors $u,v \in \mathcal{T}$ forming the basis $B = (u,v)$ of $\mathbb{E}^2$ generates a periodic lattice $\Lambda (u,v) = \{xu + yv\, |\, x,y \in \mathbb{Z}\}$ for $d_\infty$ and $\tau_{\infty}$, such that $\Lambda (u,v) \subseteq \Lambda_0$ by the definition of $\Lambda_0$. 
 For further details on translational symmetry and periodic lattice, the reader is referred, for example, to \cite[Chapters 7.3 and 9.3]{GroupSymmetry} and \cite[Chapters 1.3, 1.4, and 3.7]{TilingBook}.

\smallbreak
Two DP tangles $(\tau_\infty, \Lambda)$ and $(\tau'_\infty, \Lambda')$ (resp. two diagrams $(d_\infty, \Lambda)$ and $(d'_\infty, \Lambda')$) are said to be \textit{equal} if $\tau_\infty=\tau'_\infty$ (resp. $d_\infty = d'_\infty$) and $\Lambda=\Lambda'$ pointwise. The lattice $\Lambda (u,v)$ is isomorphic to $\mathbb{Z}^2$ and it can be generated by the standard orthonormal basis $B_0=(e_1,e_2)$ of $\mathbb{E}^2$, but also by different bases. More specifically, two bases $B=(u,v)$ and $B'=(u',v')$ generate the same point lattice $\Lambda = \Lambda (u,v) = \Lambda' (u',v')$ if and only if they are related as:
$$\begin{pmatrix}
u' \\
v' 
\end{pmatrix}
=
\begin{pmatrix}
x_1 & x_2 \\
x_3 & x_4 
\end{pmatrix}
\cdot
\begin{pmatrix}
u \\
v 
\end{pmatrix}
\textit{ , where } x_1, x_2, x_3, x_4 \in \mathbb{Z} \ \textit{ with }
(x_1 x_4 - x_2 x_3) = \pm 1.$$

The torus $T^2= \mathbb{E}^2 / \Lambda$ arises as the identification space, with respect to the boundary, of any one of the parallelograms of the periodic lattice $\Lambda$, each of which shall be called a \textit{flat torus}. Hence, a flat torus contains a tangle diagram which represents the motif diagram $d$ in $T^2$. This leads to the following definition:

\begin{definition}
  We call a motif diagram $d$ in a flat torus a \textit{flat motif diagram}, also denoted $d$. 
  Similarly, we call a motif $\tau = \tau_{\infty}$ in a thickened flat torus a \textit{flat motif}.
\end{definition}

A flat motif diagram is a particular case of a tangle diagram in the parallelogram defined by the flat torus, serving as a `tile' generating the DP diagram $d_{\infty}$ by applying the translational symmetries with respect to the basic vectors $u$ and $v$ (see Figure~\ref{DP-tangle} for an example). Clearly, flat motifs correspond bijectively to motifs.

\begin{remark}
    It is well-known that any lattice $\Lambda (u,v)$ associated to a doubly periodic structure is one of the five Bravais lattice types: square, rhombic, rectangular, hexagonal or oblique (cf. for example \cite[Chapter 9.3]{GroupSymmetry}). Each type of lattice is preserved by any isometry of the (thickened) plane. Moreover, a generating (thickened) parallelogram of each lattice type is  the identification space of a (thickened) torus. 
\end{remark}

\subsection{Minimal motifs} \label{sec:minimal-motif}

Let $d$ and $d'$ be two flat motif diagrams giving rise to the same fixed DP diagram $d_{\infty}$ and related to two distinct periodic lattices $\Lambda(u,v)$ and $\Lambda' (u',v')$, where $u'= \lambda u$ and $v'= \mu v$ with $\lambda,\mu \in \mathbb{N} \setminus \{0\}$. The flat motif diagram $d'$ is thus clearly a finite cover of the flat motif diagram $d$, and the corresponding periodic lattices satisfy the inclusion relation $\Lambda' (u',v') \subset \Lambda(u,v) \subseteq \Lambda_0$. For example, Figure~\ref{Tknot-Tlink}(c) illustrates a double cover $d'$ of the flat motif diagram $d$ in Figure~\ref{Tknot-Tlink}(b) (see also Figure~\ref{FiniteRmoves} for a six-tuple cover of a motif diagram). More generally, relations between finite covers and quotients of (flat) motif diagrams in $T^2$, carrying through to (flat) motifs in $T^2 \times I$, lead to the definition of a minimal (flat) motif (diagram) associated with a DP tangle/diagram (see Figure~\ref{crossing}(d) for an example):

\begin{definition}\label{def:minimal lattice}
A \textit{minimal (flat) motif (diagram) of} $\tau_{\infty}$ (resp. $d_{\infty}$), denoted $\tau_{min}$  (resp. $d_{min}$), is a motif which is not a finite cover of another (flat) motif (diagram) of $\tau_{\infty}$ (resp. $d_{\infty}$). The periodic lattice associated to a minimal (flat) motif diagram is a maximal periodic lattice, denoted $\Lambda_0$, and we have $d_{min}= d_{\infty} / \Lambda_0$.
\end{definition} 

In \cites{Grishanov1, Grishanov.part1, Grishanov.part2} a minimal flat motif diagram is referred to as \textit{unit cell}, while in \cites{Eleni1, Eleni2, Eleni3} it is referred to as \textit{generating cell}. A minimal motif is not unique for a given DP tangle $\tau_{\infty}$  (DP diagram $d_{\infty}$). In fact, it can be quite difficult to identify a minimal motif for $\tau_{\infty}$ not to mention a minimal motif for the whole equivalence class of $\tau_{\infty}$. The reader is referred to Subsection~\ref{sec:minimalmotif} for a further discussion on this point and an explicit example of finding a minimal motif.

\section{DP tangle isotopy} \label{sec:Dp-isotopy} 

A DP tangle $\tau_{\infty}$ is an 1-dimensional manifold embedded in the thickened plane $\mathbb{E}^2 \times I$. As such, it is allowed to undergo isotopies, that is, transformations induced by bi-continuous orientation preserving homeomorphisms  (i.e. `elastic' deformations) of $\mathbb{E}^2 \times I$. By the definition of DP tangles, we will restrict our scope only to isotopies that preserve the double periodicity, namely that they result in a DP tangle with the same topological properties.

\begin{definition} \label{def:DPequivalence}
Two DP tangles $(\tau_\infty, \Lambda)$ and $(\tau'_\infty, \Lambda')$ are called \textit{DP isotopic} if they are related via an ambient isotopy of the thickened plane taking $\tau_\infty$ to $\tau'_\infty$ and $\Lambda$ to $\Lambda'$. Two DP diagrams of two DP isotopic DP tangles shall be also  called \textit{DP isotopic}.
\end{definition}
By definition, a DP isotopy preserves the double periodicity. 

\begin{remark}\label{rem:setwise}
    We point out that a DP isotopy may result in $(\tau_\infty, \Lambda)$ and $(\tau'_\infty, \Lambda')$  being the same set-wise. This means that the two pairs can be superimposed, for example, by an automorphism of $\Lambda$. We shall expand on this point in the sequel.
\end{remark}

Isotopies that do not preserve the double periodicity of a given DP tangle could be also considered in the general topological theory of DP tangles and shall be called \textit{defect isotopies}. Yet, as we shall see below, a sequence of defect isotopies can still give rise to isotopic DP tangles. 

\begin{figure}[H]
\centerline{\includegraphics[width=5.5in]{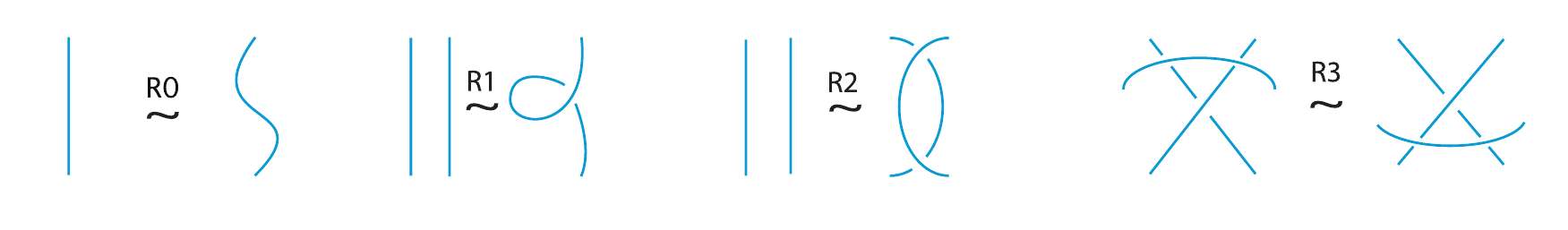}}
\caption{A planar isotopy move and the Reidemeister moves.}
\label{Rmoves}
\end{figure}

\smallbreak
Without loss of generality DP tangles are assumed to be piecewise-linear. Then, by arguments from classical knot theory, a DP tangle isotopy can be discretized by the (spatial) \textit{$\Delta$-moves}, that is, isotopies than take place within the surface of a triangle with its interior free of other arcs, and replace the arc forming one side by the other two sides or vice versa. On the diagrammatic level, these moves can be realized by local planar isotopies and Reidemeister moves, as exemplified in Figures~\ref{Rmoves} and~\ref{FiniteRmoves} (see, for example \cite{Adams}). It follows that:

\begin{proposition}\label{prop:DP-reidemeister}
    Two DP tangles $(\tau_\infty, \Lambda)$ and $(\tau'_\infty, \Lambda')$ are DP isotopic if and only if any two DP diagrams $(d_\infty, \Lambda)$ and $(d'_\infty, \Lambda')$ of theirs are related by a sequence of planar isotopies and Reidemeister moves. 
\end{proposition}

Definition~\ref{def:DPequivalence} and Proposition~\ref{prop:DP-reidemeister}, allow that a DP isotopy between two DP tangles may decompose into, or even require, intermediate isotopies which do not preserve double periodicity. 
  DP isotopies on the diagrammatic level can be classified in two categories: 
 
\begin{itemize}
    \item [-] \textit{Local DP isotopies}, generated by planar isotopies and Reidemeister moves that preserve double periodicity. This means that they can all take place simultaneously in each flat motif. In a local DP isotopy the lattice remains point-wise fixed.
    \smallbreak
    \item [-] \textit{Global DP isotopies}, corresponding to invertible, orientation-preserving affine transformations of the plane $\mathbb{E}^2$ (extending trivially along the $\mathrm{I}$-direction). These are generated by planar isotopies that may not preserve double periodicity. In a global DP isotopy the lattice is transformed to an isomorphic one.
\end{itemize}

We know that any invertible, orientation-preserving affine transformation of $\mathbb{E}^2$ is  described by a matrix of the orientation-preserving affine group \textit{Aff}$^+(2,\mathbb{R})$ of $\mathbb{E}^2$, and is generated by translations, shearings and scalings (cf. for example \cite[Chapter 1.4]{TilingBook}). These transformations include both: \textit{non-area-preserving} types, for example stretchings or contractions, and \textit{area-preserving} ones, like translations and shear deformations of the plane.   

\begin{figure}[H]
\centerline{\includegraphics[width=6in]{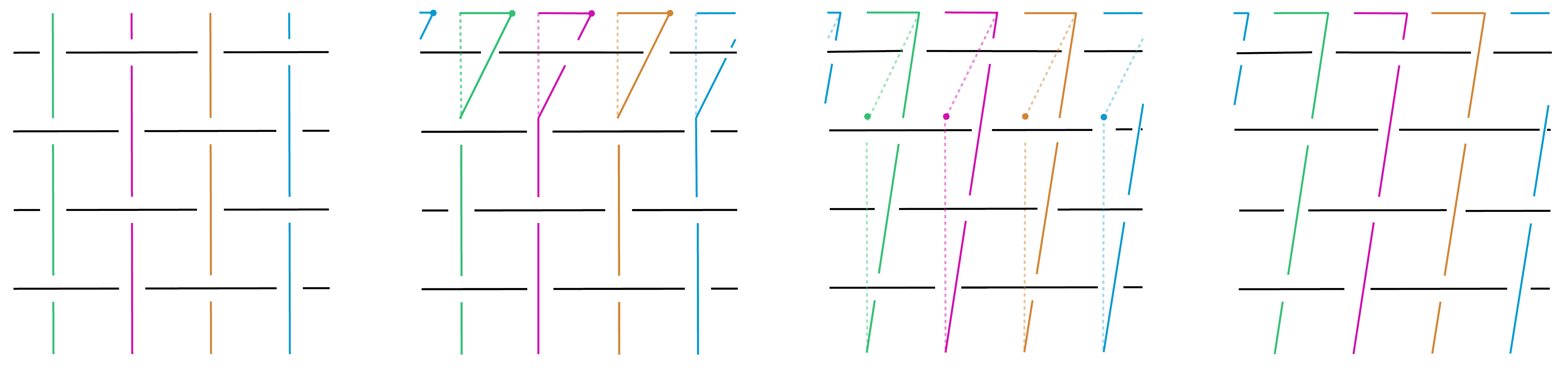}}
\caption{A shearing of a DP tangle realized by $\Delta$-moves.}
\label{Dmoves-shearing}
\end{figure}

As already pointed out, a global DP isotopy between two DP tangles $\tau_\infty$ and $\tau^\prime_\infty$  can be realized by an an infinitum of sequential planar isotopies, which in this case  may not be doubly periodic. Figure~\ref{Dmoves-shearing} illustrates an example of realizing a shearing between $\tau_\infty$ and $\tau^\prime_\infty$ by $\Delta$-moves. 

Global DP isotopies cannot be non-affine, as the periodic lattice condition would be violated. Isotopy also requires invertibility.  Moreover, they cannot be orientation-reversing as this would violate the definition of isotopy.

\begin{figure}[H]
\begin{center}
\includegraphics[width=3.5in]{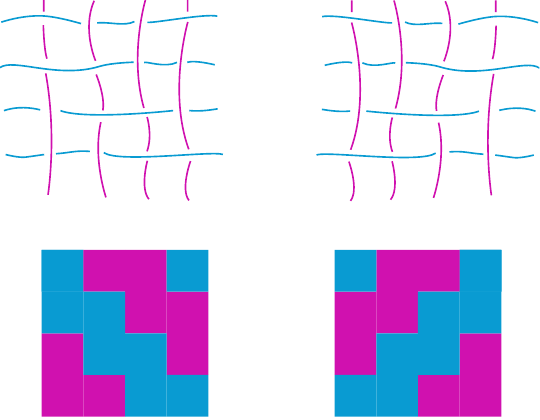}
\end{center}
\caption{Textile representation of a S-twill and Z-twill.}
\label{twill}
\end{figure}

The choice of the notion of isotopy as an equivalence relation among DP tangles is aligned with applications in materials science. In textiles for example, DP tangles are classified according to the pattern formed by their crossings at the diagrammatic level, which describes the `front side' of a material. This pattern may differ from the one on the `back side', depending on the construction method. Consider for example a particular class of DP tangles called \textit{twill weaves}, where positive (resp. negative) crossings are organized in a diagonal pattern as described in \cite{Sonia1}. If the diagonal runs in a positive slope, namely from the lower left to the upper right corner, the DP tangle is called a \textit{right-hand} twill, or \textit{Z-twill}. However, the mirror image of the DP diagram of a Z-twill gives rise to a so-called \textit{left-hand} twill, or \textit{S-twill}, where the diagonal runs in a negative slope. See  Figure~\ref{twill}. In the left-hand illustration, the blue segments (resp. pink) that for overcrossings are arranged on an S-type diagonal. On the contrary, in the left-hand illustration the over-arcs are arranged on an Z-type diagonal. These two DP tangles, which are related by an orientation-reversing transformation of the plane, are considered different in materials science. 

For the topological study of DP tangles, we would like to translate DP tangle isotopies into an equivalence relation on the level of their (flat) motifs. We will refer to this equivalence relation as \textit{DP tangle equivalence}. This is the main theme of this paper. For doing this, we shall analyze all DP local and global isotopies, taking into account all choices of lattices and bases, unlike previous studies. DP tangle equivalence has been initiated in the literature by Grishanov et al. \cite{Grishanov1} who stated a generalized Reidemeister theorem for the equivalence of DP tangles for a fixed lattice and in the context of classification of textiles. 
More recently, the third author highlighted the importance of including in the DP tangle equivalence a relation between different finite covers that lift to the same doubly periodic structure, namely the notion of \textit{scale equivalence}, for the particular case of DP isotopic weaves \cite{Sonia3}. 

\smallbreak

In the next sections we shall analyze local and global DP isotopies of DP tangles and how they translate on the level of (flat) motif diagrams. Local DP isotopies are analyzed in Section~\ref{sec:local_isotopy}, by exploiting the theory of mixed links. Global DP isotopies are addressed in Sections~\ref{sec:shiftequivalence},~\ref{sec:scaleequivalence} and~\ref{sec:dehnequivalence}, treating separately translations, scalings and shearings. In each section our strategy isto analyze the following cases: transformation  of the lattice keeping the DP diagram point-wise fixed, deformation of the DP diagram keeping the lattice point-wise fixed, which we prove to be the reverse operation,  and  simultaneous transformation of the DP diagram and its lattice, which is visible on the DP diagram level but in most cases   not visible on the level of (flat) motifs. We  also distinguish between integral and non-integral transformations, highlighting the cases where the DP diagram and/or the lattice are set-wise fixed. Moreover, in Section~\ref{sec:unimodular} we investigate further the class of unimodular linear transformations, which are all generated by shearings,  according to their classification and their effect on DP diagrams.

\section{DP tangle equivalence by local isotopies }\label{sec:local_isotopy}

In this section we present a full analysis of the local diagrammatic theory of motifs and flat motifs, taking the $(l,m)$ pair into account, for the study of equivalence of DP tangles. 

Let $\big(\tau,(l,m)\big)$ be a motif in $T^2 \times I$ that generates $\tau_{\infty}$ and $\big(d,(l,m)\big)$ a corresponding diagram of $\big(\tau,(l,m)\big)$ in the torus surface $T^2$. In the piecewise-linear setting, an isotopy of the link $\tau$ translates into a finite sequence of $\Delta$-moves in the thickened torus. On the diagrammatic level, we will see that $\Delta$-moves translate  into local moves on the diagram $d$, comprising \textit{surface isotopies} and the \textit{Reidemeister moves} (see Figure~\ref{Rmoves}), by resourcing to the theory of  mixed links.

\begin{figure}[H]
\begin{center}
\includegraphics[width=5.5in]{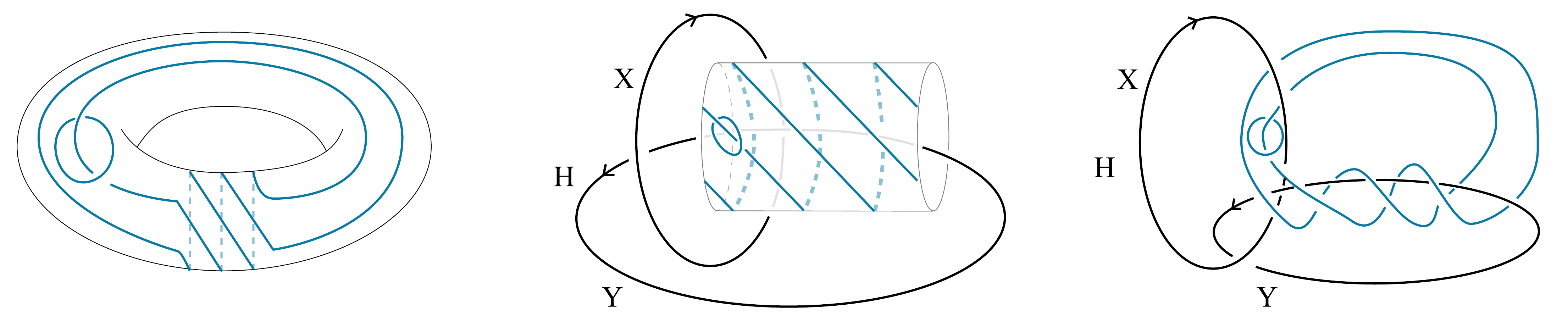}
\end{center}
\caption{A motif and its corresponding mixed link.}
\label{mixed}
\end{figure}

Let $M$ be a 3-manifold, which is the complement of a knot or link $L$ in the three-sphere $S^3$ or a closed, connected, oriented 3-manifold obtained via surgery along $L$ with additional framing. It is well-known  that isotopy classes of knots/links in $M$ correspond bijectively to isotopy classes of mixed links in $S^3$ \cite[Theorem 5.2]{LR1}. A \textit{mixed link} comprises $L$ as an oriented  point-wise \textit{fixed} sublink representing the 3-manifold $M$ and a \textit{moving} sublink representing a knot/link in $M$. In particular, a thickened torus $T^2 \times I$ can be represented as the complement in $S^3$ of the Hopf link, $H=X \cup Y$. Hence a link $\tau$ in $T^2 \times I$ corresponds in $S^3$ to the mixed link $H \cup \tau$, comprising the fixed part $H$ and the moving part representing $\tau$. For an illustration see Figure~\ref{mixed}. 

\begin{figure}[H]
\begin{center}
\includegraphics[width=5.8in]{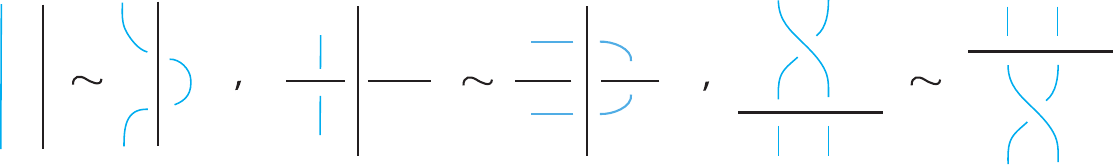}
\end{center}
\caption{\label{isom} Extended local isotopy moves for mixed links.} 
\end{figure}

On the diagrammatic level, the fixed Hopf link $H$ can be assumed to be almost flat, in the sense that arcs are coplanar and the two crossings are embedded in local 3-balls. So $H$ defines naturally a projection plane for the mixed link $H \cup \tau$ to a mixed link diagram $H \cup d$. We now recall from \cite{LR1} that  isotopy of links in the thickened torus is generated by planar isotopy, the classical Reidemeister moves for the moving part $d$, together with the extended planar isotopies that involve the fixed and the moving part of $H \cup d$, as exemplified in Figure~\ref{isom}. 

We now represent a (flat) motif $\big(\tau,(l,m)\big)$ as a mixed link $H \cup \tau$, and in analogy its corresponding (flat) motif diagram $\big(d,(l,m)\big)$ as a mixed link diagram $H \cup d$. The component $X$ of $H$ is assumed to run close and in parallel to the longitude $l$  while the component $Y$ is assumed to run close and in parallel to the meridian $m$ of $T^2$. We consider the components $X$ and $Y$  of the fixed part $H= X \cup Y$ to be ordered according to the order in the basis $B=(u,v)$ of $\mathbb{E}^2$, and oriented according to the orientation of $B=(u,v)$. 

\noindent Combining the above we conclude that:

\begin{lemma} \label{lem:mixed-link-isotopy}
Local motif isotopy classes correspond bijectively to mixed link isotopy classes.
\end{lemma}

We are now ready to translate the mixed link diagrammatic isotopy to the setting of flat motif diagrams, where $T^2$ is considered to be the identification square along $l$ and $m$,  $X$ runs close to $l$  and $Y$  runs close to $m$. In the flat picture $X$  and $Y$ form two crossed arcs. 

If a local isotopy move takes place in a mixed link, away from $X$ and $Y$, resp. in the corresponding motif diagram $\big(d,(l,m)\big)$ away from $m$ and $l$, then the same move will be visible in the interior of the flat motif diagram. 

If a local isotopy move  taking place in a mixed link involves $X$ and $Y$, the mixed link moves exemplified in Figure~\ref{isom} carry through to the corresponding motif diagram $\big(d,(l,m)\big)$. We shall examine local isotopies on the level of flat motifs, starting first with surface isotopies. These on the (flat) motif diagram level comprise:  

\begin{itemize}
    \item [a.] planar isotopies within the (flat) motif diagram,
    \smallbreak
    \item [b.] planar isotopies where an arc before lies within the motif diagram but the arc afterwards hits one boundary component (the meridian or the longitude),
    \smallbreak
    \item [c.] planar isotopies where an arc before and the arc afterwards cross one boundary component,
    \smallbreak
    \item [d.] planar isotopies where an arc before crosses one boundary component but the arc afterwards crosses both boundary components,
    \smallbreak
    \item [e.] the situation where a crossing passes through one boundary component.
\end{itemize}

\begin{figure}[H]
\centerline{\includegraphics[width=4in]{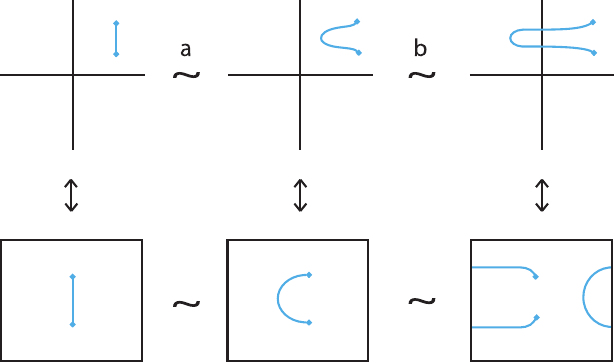}}
\vspace*{8pt}
\caption{Surface isotopies of types a) and b).}
\label{planar1}
\end{figure}

\begin{figure}[H]
\centerline{\includegraphics[width=4in]{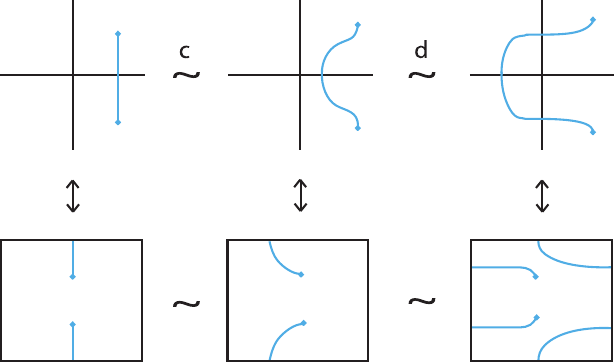}}
\vspace*{8pt}
\caption{Surface isotopies of types c) and d).}
\label{planar2}
\end{figure}

\begin{figure}[H]
\centerline{\includegraphics[width=5.7in]{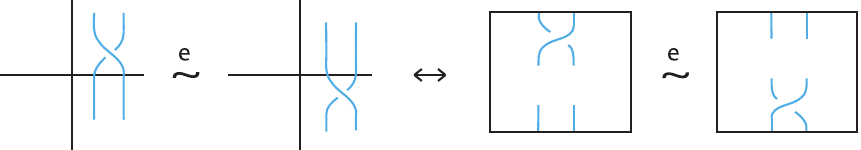}}
\vspace*{8pt}
\caption{Surface isotopies of type e).}
\label{planar3}
\end{figure}

Moves of types a) and b) are sampled in Figure~\ref{planar1}, while moves of types c) and d) are sampled in Figure~\ref{planar2}. In the figures, the top parts illustrate the moves on the motif diagram level and the bottom parts illustrate the corresponding moves on the flat motif diagram level. Moves of type e) are sampled in Figure~\ref{planar3}, where the left part shows the move on the motif diagram, while the right part shows the corresponding move on the flat motif diagram. 

\begin{figure}[H]
\centerline{\includegraphics[width=5.3in]{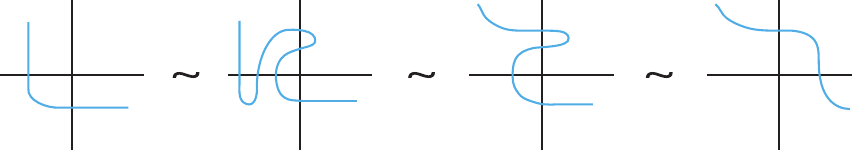}}
\vspace*{8pt}
\caption{A surface isotopy realized via the basic moves.}
\label{mR3m2b}
\end{figure}

In all figures above, the arcs may cross the  left-right or the up-down lateral sides. Also in a move type e) the crossing may be the opposite. Any other local surface isotopy can be realized by the above basic moves. As an example, the move illustrated in Figure~\ref{mR3m2b} can be realized via three moves of type b) and one move of type~d). 

\begin{figure}[H]
\centerline{\includegraphics[width=5.5in]{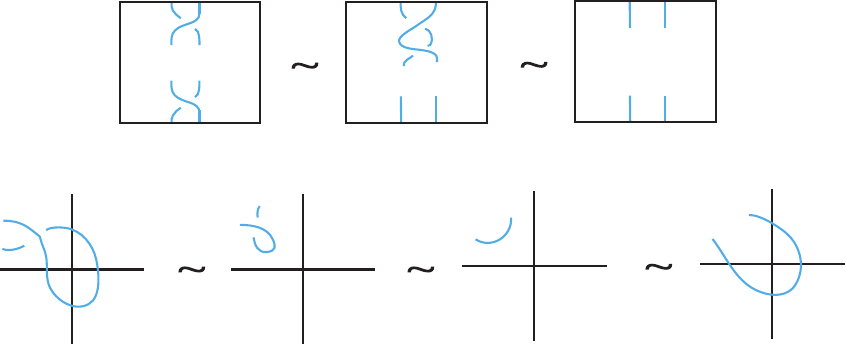}}
\vspace*{8pt}
\caption{An R2 move crossing $l$ and an R1 move crossing $m$ and $l$, both retracted to the interior of the flat motif using surface isotopies.}
\label{mixedreid}
\end{figure}

Regarding now the Reidemeister moves, we may assume that they occur far
away from $m$ and $l$ by using surface isotopies. Indeed, consider a Reidemeister move on the motif level, some arcs of which may cross $m$ or $l$ or both $m$ and $l$. Such moves are exemplified in Figure~\ref{mixedreid}. Using surface isotopies, both sides of the Reidemeister move can be pushed away from $m$ and $l$, which results in the move to take place in the interior of the flat motif  diagram. After the move is performed we push back all arcs to their original position using surface isotopies.

The above exhaustive analysis, in view of \cite[Theorem 5.2]{LR1} and Lemma~\ref{lem:mixed-link-isotopy}, establishes the following:

\begin{lemma}\label{lem:localmotifisotopy}
Two (flat) motifs are isotopic if and only if two (flat) motif diagrams of theirs differ by a finite sequence of surface isotopy moves and the classical Reidemeister moves, as exemplified in Figures~\ref{Rmoves},~\ref{planar1},~\ref{planar2}, and ~\ref{planar3}. 
\end{lemma}

In  motif isotopy we have a moving link in the thickened torus, which is assumed fixed. So, the moves in Lemma~\ref{lem:localmotifisotopy} on a motif diagram $d$ generate, in turn, on the corresponding DP diagram $d_{\infty}$ (local) planar isotopies and Reidemeister moves, that preserve the double periodicity,  with the same fixed lattice $\Lambda(u,v)$. These are the \textit{local DP planar isotopies} and the \textit{DP Reidemeister moves}. The above lead now to the following.

\begin{proposition}\label{prop:localmotifisotopy}
Two DP tangles are locally isotopic if and only if two (flat) motifs of theirs are isotopic. Equivalently, if and only if  they differ by a finite sequence of local DP planar isotopies and DP Reidemeister moves.
\end{proposition}

\begin{figure}[H]
\centerline{\includegraphics[width=3.5in]{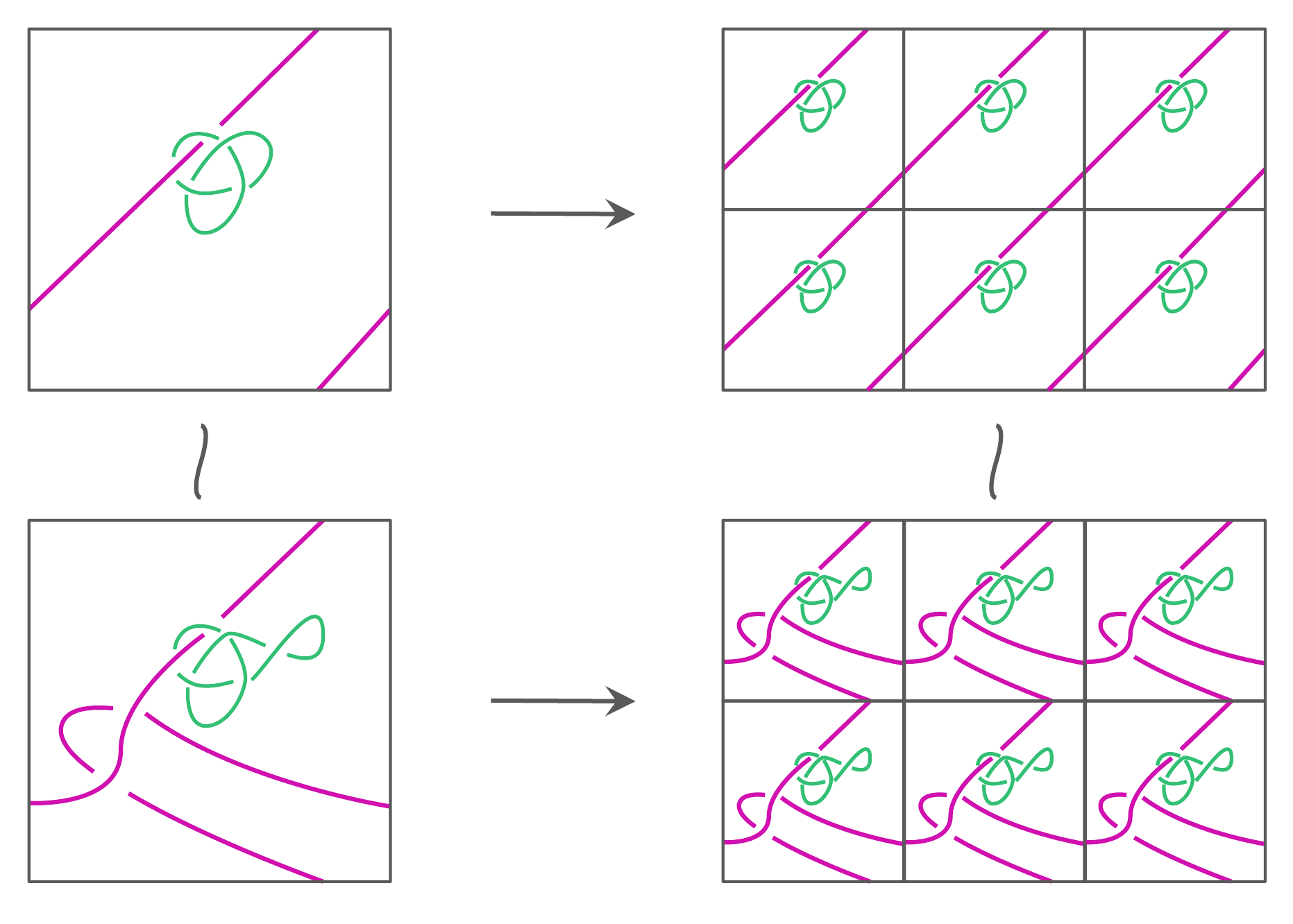}}
\caption{\label{FiniteRmoves} On the top, a flat motif and a corresponding (finite) cover. On the bottom, a locally isotopic flat motif and the corresponding DP Reidemeister moves.}
\end{figure}

\begin{remark}
Proposition~\ref{prop:localmotifisotopy} applies to any finite cover of the corresponding (flat) motif $\tau$. An example is illustrated in Figure~\ref{FiniteRmoves}. As the figure suggests, if we take the same finite cover of a given motif before and after performing a local move, then this yields isotopy on the level of the two covers, since they differ by a finite sequence of that repeated  local move.
\end{remark}

Concluding, isotopy moves between (flat) motif diagrams on the fixed torus $T^2$ correspond to local isotopy moves between DP diagrams supported by the same lattice,  meaning that the longitude-meridian pair $(l,m)$  remains fixed.  We may also consider moving the $(l,m)$ pair keeping the link fixed or, more generally, transforming the lattice. For example, stretches and contractions of flat motif diagrams, or in other words, torus inflations and deflations, are not included in the local isotopy discussion, as these also transform the underlying lattice. These transformations, along with all global isotopies are addressed in the sections below. Since global isotopies result from affine transformations of the plane $\mathbb{E}^2$, which are assumed to extend to $\mathbb{E}^2 \times I$ \, trivially along the $\mathrm{I}$-direction, we will focus our analysis exclusively on the diagrammatic level.

\section{DP tangle global equivalence: translation and shift equivalence}\label{sec:shiftequivalence}

Let $(d_\infty,\Lambda)$  be a DP diagram of a DP tangle $(\tau_\infty,\Lambda)$ with its periodic lattice. Let also $d$ denote the (flat) motif diagram of $d_{\infty}$ in the torus $T^2$ equipped with the longitude-meridian pair $(l,m)$. In this section we examine the effect of translation in the context of DP tangle global isotopy. A \textit{translation} of $\mathbb{E}^2$ is an affine transformation which is an isometry that shifts each point of the plane by any fixed vector $t$.  A \textit{translation of the lattice} $\Lambda$ is the restriction of a plane translation on the points of $\Lambda$. Clearly, a lattice translation preserves the Bravais lattice type. A \textit{translation of the DP diagram} $d_\infty$ is the restriction of a plane translation on the points of $d_\infty$. We shall analyze below the different cases of translations, paying special attention to the distinction between integral and non-integral ones.

\bigbreak
\noindent \textit{Case 1: Translation of the lattice keeping the DP diagram point-wise fixed.} 
\smallbreak

\noindent A translation of the lattice $\Lambda$ by a fixed vector $t$ results in a different integer lattice $\Lambda'$. Keeping the DP diagram $d_\infty$ point-wise fixed, the lattice translation gives rise to a new DP diagram $(d_\infty,\Lambda')$ and thus a new DP tangle $(\tau_\infty,\Lambda')$. Indeed, the coordinate system $(\mathrm{O},u,v)$ can  be shifted by any vector $t$  in the plane, while keeping $d_\infty$ fixed, and still preserve the translational symmetry of the $d_\infty$ by definition of a DP tangle, so the shifted lattice $\Lambda (u+t,v+t) = \{x(u+t) + y(v+t)\, |\, x,y \in \mathbb{Z}\}$ is still a lattice of periodicity for $d_\infty$.

\begin{figure}[H]
\centerline{\includegraphics[width=4in]{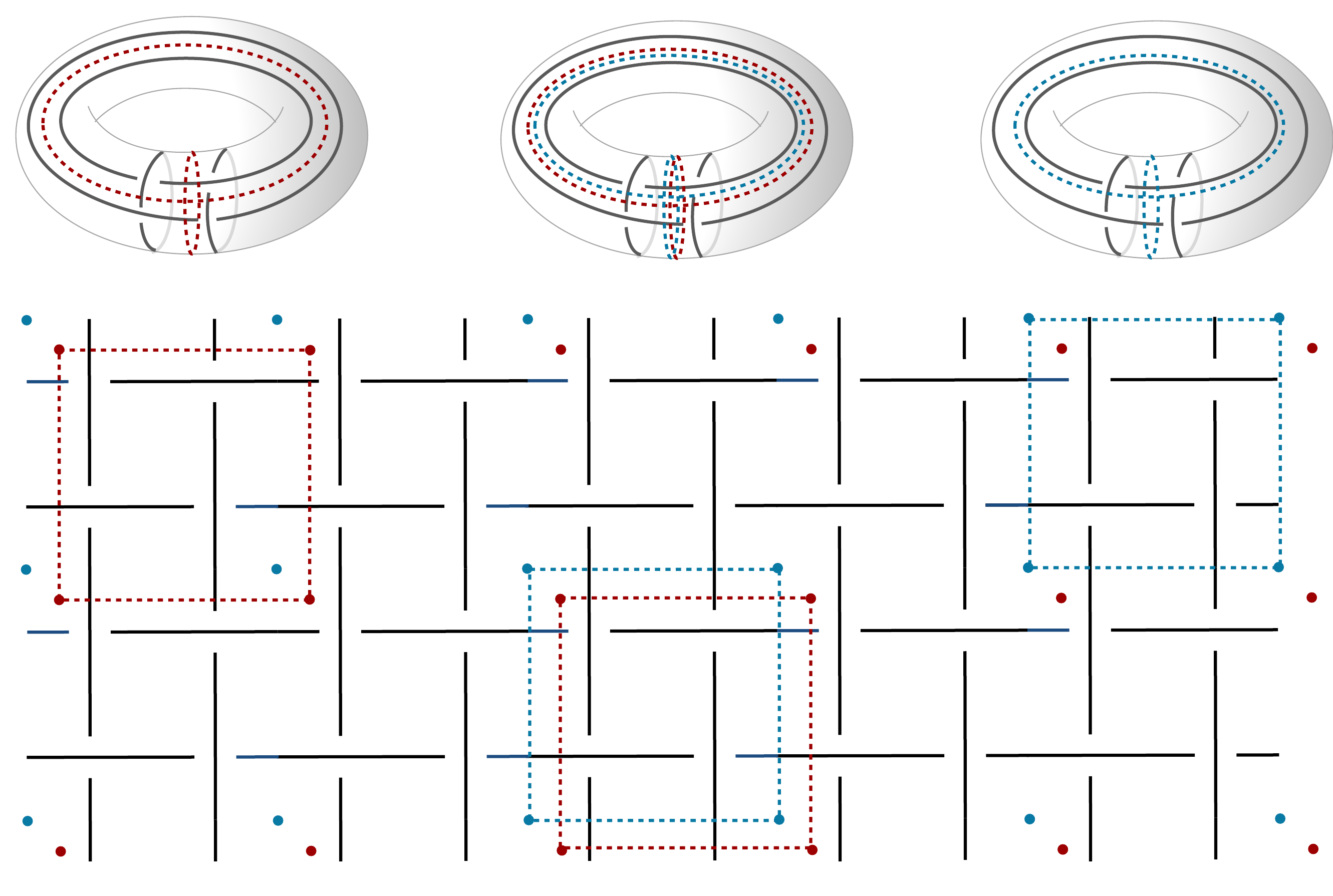}}
\caption{\label{fig:Shift-eq} On the top, the motif diagrams $(d,(l,m))$ (left), $(d,(l',m'))$ (right), and their superposition (middle). On the bottom, their corresponding flat motifs in the DP diagram $d_\infty$ with the corresponding red and blue lattices.}
\end{figure}

On the level of the motif diagram this operation clearly corresponds to a shift of $(l,m)$ to a new longitude-meridian pair $(l',m')$. Hence, it generates two  flat motif diagrams $d$ and $d'$, corresponding to the motif diagrams $(d,(l,m))$ and $(d,(l',m'))$,  related via a translation by the vector $t$ in $d_\infty$, as illustrated in Figure~\ref{fig:Shift-eq}. In the figure, we can see a longitude-meridian pair in either side, and their superposition in the middle illustration, highlighting the translation. In particular, considering the longitude-meridian pairs for $d$ on $T^2$, $(l,m)$, $(l',m)$ and $(l,m')$ we observe that: the pair $(l',m)$ indicates a vertical shift of the (flat) motif diagram $d$ created by the pair $(l,m)$, while the pair $(l,m')$ indicates a horizontal shift of the (flat) motif diagram $d$. Any shift $(l',m')$ can result as the combination of a horizontal and a vertical shift, where the order is not relevant. 

We note that shifts of longitude-meridian pairs include also shifts by multiples of $2\pi$. These correspond to \textit{integral translations} or \textit{$\mathbb{Z}$-translations}  of the underlying lattice $\Lambda$, fixing $\Lambda$ set-wise. On the motif level, an integral translation of $\Lambda$ gives rise to the same point-wise (flat) motif diagram $(d,(l,m))$. Hence, an integral translation is invisible on the flat motif level. On the other hand, a \textit{ non-integral translation} or \textit{ non-$\mathbb{Z}$ translation} of $\Lambda$ gives rise to different longitude-meridian pairs, and thus different motif diagrams $(d,(l,m))$ and $(d,(l',m'))$ and different corresponding flat motif diagrams $d$ and $d'$. 

\bigbreak
\noindent \textit{Case 2: Translation of the DP diagram  keeping the lattice point-wise fixed.} 
\smallbreak

\noindent A translation of the DP diagram $(d_\infty,\Lambda)$ is a displacement of each point by a fixed vector $t$. Keeping the lattice $\Lambda$ point-wise fixed, the translation gives rise to a new DP diagram $(d'_\infty,\Lambda)$ and thus a new DP tangle $(\tau'_\infty,\Lambda)$. This case is the reverse operation of Case 1. The result can be equivalently achieved  by fixing the DP diagram $d_\infty$ and translating the lattice $\Lambda$ by the fixed vector $-t$.

\begin{figure}[H]
\centerline{\includegraphics[width=1.5in]{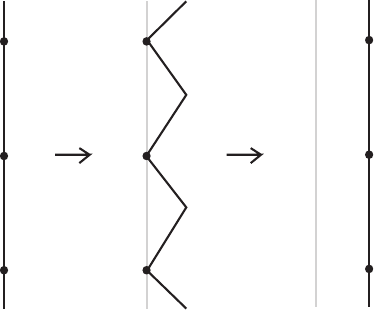}}
\vspace*{8pt}
\caption{\label{translation-Delta} Attaining DP diagram translation by sequential isotopy moves.}
\end{figure}

What is interesting in this case is that the DP diagram translation can be attained  by an infinitum of sequential planar isotopies, which in this case can be doubly periodic.  On the level of the motif diagrams the translation of $d_\infty$  corresponds to transforming the motif diagram $(d,(l,m))$ to a new motif diagram $(d',(l,m))$  via a finite sequence of surface isotopy moves. On the level of flat motif diagrams, such moves may involve some arcs or crossings intersecting $l$ or $m$ or both $l$ and $m$. By the analysis in Section~\ref{sec:local_isotopy} these are the surface isotopy moves (a)--(e), recall Figure~\ref{planar1}, ~\ref{planar2}, ~\ref{planar3}. An abstraction of the diagram translation is shown in Figure~\ref{translation-Delta}, where an arc is shifted to a parallel arc by $\Delta$-moves. 

In the special case of a $\mathbb{Z}$-translation of the DP tangle, this  fixes the DP tangle set-wise, keeping the lattice point-wise fixed.  A $\mathbb{Z}$-translation is invisible on the motif and the flat motif level.  

\bigbreak
\noindent \textit{Case 3: Simultaneous translation of the DP diagram and the lattice.} 
\smallbreak

\noindent A simultaneous translation of the DP diagram and its lattice $(d_\infty,\Lambda)$ is visible on the level of the DP tangle, but not visible on the level of (flat) motif diagrams. It results in a new DP diagram and lattice $(d'_\infty,\Lambda')$ both translated by the same vector $t$.  
In the special case of $\mathbb{Z}$-translation, the pair $(d_\infty,\Lambda)$ is set-wise fixed, in the sense that each one the pairs $d_\infty, d'_\infty$ and $\Lambda,\Lambda'$ can be  superimposed, while for a non-$\mathbb{Z}$-translation they cannot be superimposed.

\smallbreak

The above  lead to the following definition.

\begin{definition}\label{def:shift-equivalent}
Two (flat) motif diagrams of a DP tangle are said to be \textit{shift equivalent} if they are related by a shift of the $(l,m)$ pair, equivalently, by a specific finite sequence of surface isotopy moves.
\end{definition} 

So, by the  analysis above and Definition~\ref{def:shift-equivalent} we derive the following:

\begin{proposition}\label{prop:shift-equivalent}
Two DP tangles are related by a DP translation if and only if two (flat) motif diagrams of their are shift equivalent.
\end{proposition} 

\begin{remark} 
In view of shift equivalence, we can conclude that any surface isotopy move of the motif diagram that crosses the longitude or/and the meridian $(l,m)$ can be realized away from $(l,m)$, that is, in the interior of the flat motif diagram. Indeed, since the isotopy move is local, by an appropriate lattice translation of Case 1, equivalently, by an appropriate shift of $(l,m)$, the isotopy move can occur away from $(l,m)$. Doing the opposite shift/the  opposite lattice translation after realizing the move, brings $(l,m)$/the lattice to its initial position. 
\end{remark}

\section{DP tangle equivalence: scaling and scale equivalence}\label{sec:scaleequivalence}

Let  $\mathbb{E}^2$ be the Euclidean plane equipped with a basis $B=(u,v)$. Let further $(d_\infty,\Lambda)$  be a DP diagram of a DP tangle $(\tau_\infty,\Lambda)$ with its periodic lattice. Let also $d$ denote the (flat) motif diagram of $d_{\infty}$ in the torus $T^2$ equipped with the longitude-meridian pair $(l,m)$. In this section we examine the effect of scaling in the context of DP tangle global isotopy. A \textit{scaling} of $\mathbb{E}^2$ is a linear transformation that stretches or contracts the plane by a scale factor in each direction of the basis, thus preserving the angle between the two basis vectors. 

A \textit{scaling of the lattice} $\Lambda(u,v)$  is the restriction of a plane scaling on the points of $\Lambda$. This operation sends the basis $(u,v)$ of $\mathbb{E}^2$ to a re-scaled basis $(u',v') = (\lambda u, \mu v)$, where the scale factors $\lambda, \mu$ are fixed positive real numbers. A scaling of the lattice $\Lambda(u,v)$ results in a different integer lattice $\Lambda'(u',v')$ which may not be of the same Bravais lattice type. In particular, we observe that the pair $(\lambda u, v)$ indicates a horizontal scaling of the lattice, while the pair $(u, \mu v)$ indicates a vertical scaling of the lattice. Then, any scaling can result as the combination of a horizontal and a vertical one, where the order is not relevant. 

\begin{figure}[H]
\centerline{\includegraphics[width=4in]{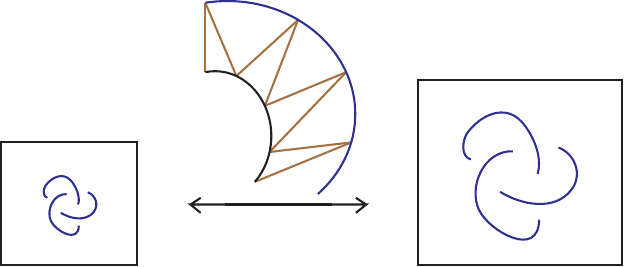}}
\vspace*{8pt}
\caption{\label{scaling-by-Delta} Attaining DP diagram scaling by sequential isotopy moves.}
\end{figure}

A \textit{scaling of the DP diagram} $(d_\infty,\Lambda)$  is the restriction of a plane scaling (stretching/ contraction) on the subset $d_\infty$. A linear transformations preserves double periodicity, so this operation gives rise to a new  diagram $d'_\infty$ which is also double periodic with a lattice $\Lambda'$. A DP diagram scaling can be attained  by an infinitum of sequential isotopies, as suggested in Figure~\ref{scaling-by-Delta},  which in this case will, in general, not be doubly periodic.  

We shall analyze below the different cases of scaling, especially with regard to the distinction between integral and non-integral ones.

\bigbreak
\noindent \textit{Case 1: Scaling of the lattice keeping the DP diagram point-wise fixed.} 
\smallbreak

\noindent  We shall first consider $\mathbb{N}$-scalings, namely scalings where the basis $(u,v)$ of $\mathbb{E}^2$ is sent to a re-scaled basis $(u',v') = (\lambda u, \mu v)$, where $\lambda, \mu \in \mathbb{N}$ are fixed positive integers. Keeping the DP diagram $d_\infty$ point-wise fixed, the integral lattice scaling gives rise to a new DP diagram $(d_\infty,\Lambda')$ and thus a new DP tangle $(\tau_\infty,\Lambda')$. Indeed, the coordinate system $(\mathrm{O},\lambda u,\mu v)$ defines new translational symmetry for the fixed $d_\infty$ along the vectors $u'=\lambda u$ and $v'=\mu v$, giving rise to the periodic lattice $\Lambda'(u',v')=\Lambda'(\lambda u,\mu v) = \{x\lambda u + y\mu v\, |\, x,y \in \mathbb{Z}\}$, being a sub-lattice of $\Lambda (u,v) = \{xu + yv\, |\, x,y \in \mathbb{Z}\}$.

\begin{figure}[H]
\centerline{\includegraphics[width=5.5in]{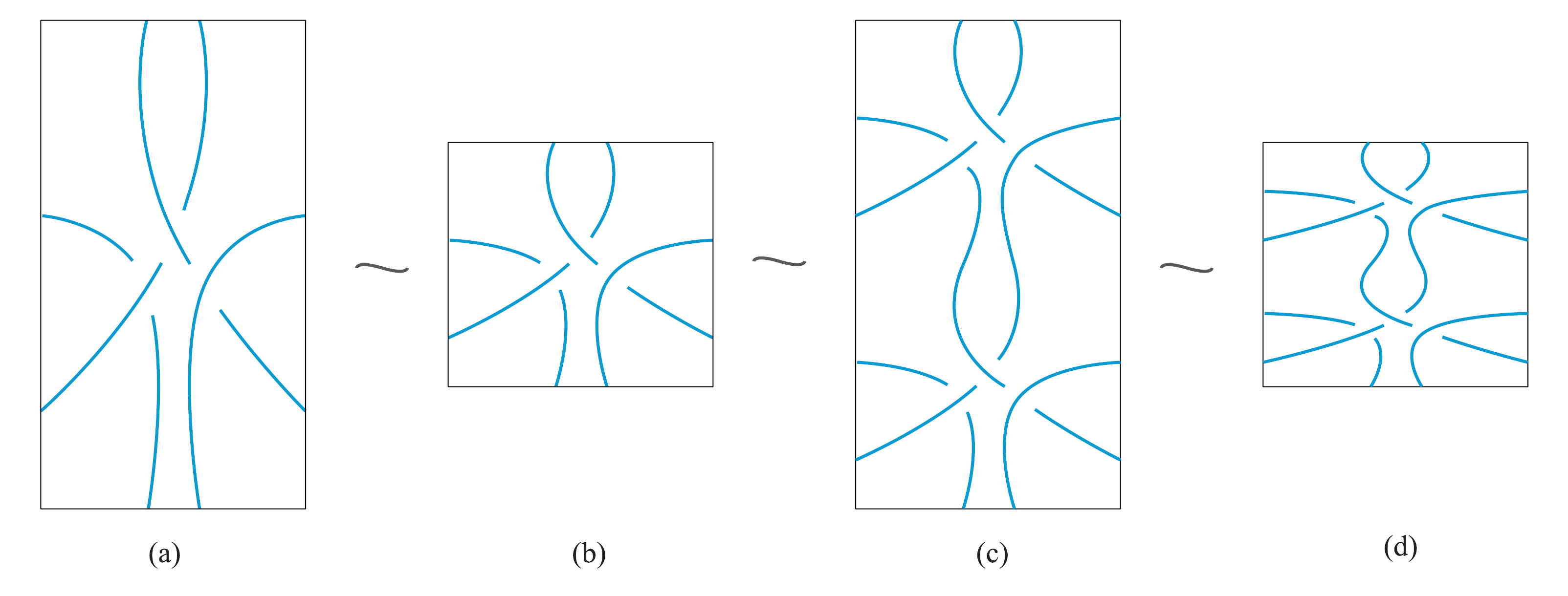}}
\vspace*{8pt}
\caption{\label{Tknot-Tlink} Motif (a) is a stretching of (b);  motif (d) is a contraction of (c); motifs (b) and (c) are scale equivalent.}
\end{figure}

On the level of the motif diagram, this operation corresponds to a finite covering map of the torus $T^2$: $\bar{\rho}: \overline{T^2} \rightarrow{} T^2$ sending the longitude-meridian pair $(l,m)$ of the torus $T^2$ to the longitude-meridian pair $(l',m')$ of the torus $\overline{T^2}$. In particular, the tori $T^2$ and $\overline{T^2}$ are related to the two bases via the covering maps $\rho: \mathbb{E}^2 \rightarrow{} T^2$ assigning the longitude $l$ of $T^2$ to $u$ and the meridian $m$ of $T^2$ to $v$, and $\rho': \mathbb{E}^2 \rightarrow{} \overline{T^2}$ assigning the longitude $l'$ of $\overline{T^2}$ to $u'$ and the meridian $m'$ of $\overline{T^2}$ to $v'$. The covering map $\bar{\rho}$ is said to be an \textit{angle-preserving covering map}. Hence, this operation generates two different flat motif diagrams $d$ and $d'$, corresponding to the motif diagrams $(d,(l,m))$ and $(d',(l',m'))$.  Figure~\ref{Tknot-Tlink}(b) and (c) illustrates an example of vertical lattice stretching by a factor 2 from (b) to (c). In the figure, we can see a flat motif diagram (b) and an angle-preserving double cover (c) of (b). 

A non-$\mathbb{N}$ scaling of the lattice $\Lambda(u,v)$ to a lattice $\Lambda'(u',v')$, where  $(u',v') = (\lambda u, \mu v)$, for scale factors $\lambda, \mu$ arbitrary positive, non-integer real numbers in the open interval $(0,1)$,  while keeping the DP diagram $d_\infty$ point-wise fixed,  will not preserve in general the double periodicity, since $u'$ and $v'$ will not, in general, be vectors of translational symmetry for the fixed $d_\infty$. (The general case is covered by modularity.) However, there are cases where scalings by rational numbers will work. For example, in Figure~\ref{Tknot-Tlink}(b) and (c) the scaling from (c) to (b) is a vertical contraction of the lattice by a factor of $\frac{1}{2}$. On the motif level, (c) is a double cover of (b), equivalently, (b) is a quotient of (c). Even in this case, the scaling from (c) to (b) works because (c) is not a minimal motif diagram.

\begin{remark} 
A scaling by a rational number cannot occur in the case of a minimal (flat) motif diagram, since, by Definition~\ref{def:minimal lattice}, it cannot be divided into a `smaller' (flat) motif diagram for the same DP diagram. 
\end{remark}

In conclusion, only $\mathbb{N}$-scalings of the lattice $\Lambda(u,v)$ will preserve double periodicity of the fixed DP diagram $d_\infty$. 
The above lead to the following definition, generalizing the notion of \textit{scale equivalence} of DP weaves from \cite{Sonia1}.

\begin{definition}\label{def:scale-equivalent}
Two (flat) motif diagrams of a DP tangle are said to be \textit{scale equivalent} if they are related by an angle-preserving covering map.
\end{definition}

\bigbreak
\noindent \textit{Case 2: Simultaneous scaling of the DP diagram and the lattice.} 
\smallbreak

\noindent A simultaneous scaling of the DP diagram and its lattice $(d_\infty,\Lambda)$ results from a scaling of $\mathbb{E}^2$ by scale factors $\lambda, \mu \in \mathbb{R}_+$ positive real numbers in each direction of the basis, applied to both $\Lambda(u,v)$ and  $d_\infty$, that is, carrying along the DP diagram and its lattice. The resulting pair is a new DP diagram $(d'_\infty,\Lambda')$ and thus corresponds to a new DP tangle $(\tau'_\infty,\Lambda')$.

This operation is  visible on the level of the DP diagram, resulting in a corresponding scaling of the flat motif. On the level of  motif diagrams, it is topologically trivial, since it corresponds to a torus inflation or deflation, which is a global isotopy of the torus, carrying along the link diagram. This generates two different motif diagrams $(d,(l,m))$ and $(d',(l',m'))$ corresponding to two different flat motif diagrams $d$ and $d'$. Two examples are illustrated in Figure~\ref{Tknot-Tlink}: the flat motif diagram (a) is a stretching of the flat motif diagram (b) generated by torus inflation, while conversely (b) is a contraction of (a) induced by torus deflation. An analogous example is between the flat motif diagrams (c) and (d). 

Two flat motif diagrams are said to be related by a \textit{stretching/contraction} if their corresponding motif diagrams are related by a \textit{torus inflation/deflation}. In conclusion, any stretching or contraction of a flat motif diagram, along with the lattice, gives rise to isotopic DP diagrams, and thus isotopic DP tangles.

\bigbreak
\noindent \textit{Case 3: Scaling of the DP diagram keeping the lattice point-wise fixed.} 
\smallbreak

\noindent Consider first a contraction of $d_\infty$ by scale factors $\frac{1}{\lambda}, \frac{1}{\mu}$ along the two basis vector directions respectively, for $\lambda, \mu \in \mathbb{N}$ positive integers. Keeping the lattice $\Lambda$ point-wise fixed, this operation gives rise to a new DP diagram $(d'_\infty,\Lambda)$ and thus a new DP tangle $(\tau'_\infty,\Lambda)$. 

On the level of flat motifs, this operation corresponds to a contraction of the flat motif $d$ by scale factors $\frac{1}{\lambda}, \frac{1}{\mu}$. This corresponds to a covering map  $\bar{\rho'}: T^2 \rightarrow{} T^2$ of the torus $T^2$ sending the longitude-meridian pair $(l,m)$ of the torus to the same longitude-meridian pair $(l,m)$ of the same torus $T^2$ since the lattice and the basis of $\mathbb{E}^2$ are fixed. Such a covering map is also angle-preserving covering. Hence, this operation generates two different flat motif diagrams $d$ and $d'$, corresponding to the motif diagrams $(d,(l,m))$ and $(d',(l,m))$. Examples are illustrated in Figure~\ref{Tknot-Tlink}(b) and (d) for a fixed lattice, as well as (a) and (c) for another example of fixed lattice. In the figure, we can see a flat motif diagram (b) and an angle-preserving double cover (d) of (b) for a fixed lattice. From the viewpoint of DP diagram scaling, we have a contraction of the DP diagram from (b) to (d), and conversely a stretching of the DP diagram from (d) to (b). An analogous example is between the flat motif diagrams (a) and (c). 

Up to Case 2, Case 3 can be seen as the reverse operation of Case 1: the result of contracting a DP diagram by scale factors $\frac{1}{\lambda}, \frac{1}{\mu}$, keeping the lattice fixed, can be equivalently achieved by a combination of Cases 1 and 2, namely, by first scaling the lattice by the inverse factors $\lambda, \mu$, keeping the DP diagram fixed (Case 1), and then scaling simultaneously the DP diagram and the lattice by the initial scale factors $\frac{1}{\lambda}, \frac{1}{\mu}$ (Case 2). See, for example, Figure~\ref{Tknot-Tlink} on the level of flat motifs.  
From (b) to (d) we have a DP contraction by the factors $1, \frac{1}{2}$ while keeping the lattice fixed. The result can be realized equivalently by first passing  from (b) to (c) by a lattice  scaling by the inverse factors $1, 2$ keeping the DP diagram fixed (Case 1), and then going from (c) to (d) by a simultaneous scaling by the scale factors $1, \frac{1}{2}$ (Case 2).  
 In conclusion, this case  also demonstrates scale equivalence as in Definition~\ref{def:scale-equivalent} in Case 1. 

An arbitrary scaling, that is, a stretching or contraction, of the DP diagram $(d_\infty,\Lambda)$ keeping the lattice $\Lambda(u,v)$ point-wise fixed, and thus the basis fixed, will not, in general, give rise to a DP diagram $(d'_\infty,\Lambda)$. Even a $\mathbb{N}$-scaling of the DP diagram $(d_\infty,\Lambda)$, that is, an integral stretching or contraction of $d_\infty$ by positive integers, keeping the lattice $\Lambda(u,v)$ point-wise fixed, will in general not result in a DP diagram, unless the corresponding flat motif diagram possesses internal translational symmetries. View, for example, the flat motif diagram (c) in Figure~\ref{Tknot-Tlink}. A vertical stretching of the flat motif diagram (c) by a factor 2, keeping the lattice fixed, will give rise to the flat motif diagram (a),  which is a flat motif diagram for an isotopic DP diagram $d'_\infty$ arising by stretching vertically by the factor 2 the initial DP diagram $d_\infty$. 

\begin{remark} 
 In the particular  case where our starting flat motif diagram $d$ is a minimal motif, so $\Lambda$ is a maximal lattice, then any  $\mathbb{N}$-stretching of $d$, keeping $\Lambda$ fixed, will not correspond to a DP diagram, by the definition of a minimal motif.
\end{remark} 

By the analysis of all three cases above and Definition~\ref{def:scale-equivalent} we derive the following:

\begin{proposition}\label{prop:scale-equivalent}
Two DP tangles are related by a scaling  if and only if two (flat) motif diagrams of their are scale equivalent and/or are related by torus inflation/deflation.
\end{proposition} 

\begin{remark}
  The notion of scale equivalence implies that a motif containing a single component, i.e. a knot, may be scale equivalent to a motif with multiple components, that is, a link. As an example, consider Figure~\ref{Tknot-Tlink}. On the left (instances (a) and (b)), a motif that corresponds to a knot in $T^2 \times I$ is illustrated, while on the right (instances (c) and (d)), a motif of the same DP tangle is presented, which now corresponds to a link in $T^2 \times I$.
\end{remark}

\section{DP tangle equivalence: shearing and Dehn equivalence}\label{sec:dehnequivalence}

As usual,  $\mathbb{E}^2$ denotes the Euclidean plane equipped with a basis $B=(u,v)$. Further $(d_\infty,\Lambda)$ denotes a DP diagram of a DP tangle $(\tau_\infty,\Lambda)$ with its periodic lattice and $d$  the (flat) motif diagram of $d_{\infty}$ in the torus $T^2$ equipped with the longitude-meridian pair $(l,m)$. In this section we examine the effect of shearing in the context of DP tangle global isotopy. A \textit{shearing} of $\mathbb{E}^2$ is an area-preserving linear transformation that shears the plane by a shear factor in a single direction of the basis, thus changes the angle between the two basic vectors. More precisely, a shearing sends the basis $B= (u,v)$ of $\mathbb{E}^2$ to a sheared basis $B'=(u',v')$ related by a specific type of matrix in $SL(2,\mathbb{R})$. Namely, the matrices $M_k$ below indicate a horizontal shearing of the basis, while the matrices $L_k$ indicate a vertical shearing of the basis, where $k\in \mathbb{R}$: 
\begin{equation}
\begin{aligned}
M_k &= \begin{pmatrix}
1 & k \\
0 & 1
\end{pmatrix}, &
L_k &= \begin{pmatrix}
1 & 0 \\
k & 1
\end{pmatrix}.
\end{aligned}
\label{shearing-matrices}
\end{equation} 
A \textit{shearing of the lattice} $\Lambda(u,v)$ is the restriction of a plane shearing to the points  of $\Lambda$. We will distinguish between integral and non-integral shearings. An integral shearing, or $\mathbb{Z}$-shearing, of the lattice $\Lambda(u,v)$ is realized via matrices $M_k, \, L_k \in SL(2,\mathbb{Z})$, for $k\in \mathbb{Z}$.  This results in  the same integer point lattice but with a different basis, $\Lambda(u',v')$. Namely, $\Lambda' = \Lambda(u',v') = \Lambda(u,v) = \Lambda$ set-wise. A non-integral shearing, or non $\mathbb{Z}$-shearing, does not have this nice property and is discussed further below.

A \textit{shearing of the DP diagram} $(d_\infty,\Lambda)$ is the restriction of a plane shearing on the subset $d_\infty$. A linear transformation  preserves double periodicity, so this operation gives rise to a new  diagram $d'_\infty$ which is also doubly periodic with a possibly different lattice $\Lambda'$. A shearing of a DP diagram can be attained by an infinitum of sequential isotopies, as suggested in Figure~\ref{Dmoves-shearing}, which in this case will, in general, not preserve double periodicity.  

We shall analyze below the different cases of shearing, especially with regard to the distinction between integral and non-integral ones.

\bigbreak

\noindent \textit{Case 1: Shearing of the basis keeping the DP diagram point-wise fixed.} 
\smallbreak

\noindent  We shall first consider $\mathbb{Z}$-shearings of the lattice, namely shearings where the basis $B=(u,v)$ of $\mathbb{E}^2$ is sent to a sheared basis $B' = (u',v') = A (u,v)$, where $A$ is a matrix of type $M_k$ or $L_k$, with $k\in \mathbb{Z}$, as defined in~(\ref{shearing-matrices}). Keeping the DP diagram $d_\infty$ point-wise fixed, the change of basis gives rise to the same DP diagram $(d_\infty,\Lambda)$ and thus the same DP tangle $(\tau_\infty,\Lambda)$. So, the change of basis is invisible on the level of DP tangle and point lattice.

\begin{figure}[H]
\centerline{\includegraphics[width=6in]{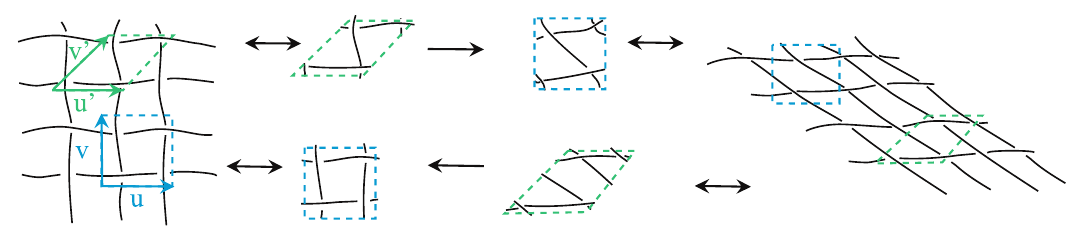}}
\vspace*{8pt}
\caption{\label{shearing} 
A DP diagram with a fixed lattice and two different bases of $\mathbb{E}^2$, leading to a shearing of the DP tangle.}
\end{figure}

On the level of the flat motif diagrams, this operation changes the geometry of the parallelogram flat torus sending the longitude-meridian pair $(l,m)$ of the torus $T^2$ to the longitude-meridian pair $(l',m)$ or $(l,m')$ of $T^2$, as illustrated in the very left side of Figure~\ref{shearing}. The figure illustrates an example of vertical lattice shearing from the blue basis  to the green basis, and the corresponding flat motif diagrams. Hence, this operation generates two different flat motif diagrams $d$ and $d'$, corresponding to the motif diagrams $(d,(l,m))$ and $(d',(l',m))$ or $(d',(l,m'))$. See second illustrations of Figure~\ref{shearing}. 

We shall next consider non-integral, or non-$\mathbb{Z}$, shearings of the basis vectors $u,v$ and consequently of the lattice $\Lambda(u,v)= \{xu + yv\, |\, x,y \in \mathbb{Z}\}$,  where the basis $B=(u,v)$ of $\mathbb{E}^2$ is sent to a sheared basis $B' = (u',v') = A (u,v)$, where $A$ is a matrix of type $M_k$ or $L_k$, with $k\in \mathbb{R}\setminus \mathbb{Z}$, as defined in~(\ref{shearing-matrices}). To preserve the double periodicity, the image $\Lambda'(u',v')$ of $\Lambda(u,v)$  must still be a periodic lattice for $d_\infty$. Thus, we need to find conditions on the sheared basis vectors $u',v'$ so that they are elements of the translation group $\mathcal{T}$ of $d_\infty$, or in other words so that $\Lambda' \subseteq \Lambda_0$, where $\Lambda_0 (u_0,v_0) = \{xu_0 + yv_0\, |\, x,y \in \mathbb{Z}\}$ is the maximal periodic lattice of $d_\infty$ (see Section~\ref{sec:lattice}).

\smallbreak

We consider the case of a shearing of type $M_k$ applied to $\Lambda (u,v)$. Since $\Lambda \subseteq \Lambda_0$, the vectors $u$ and $v$ can be written in terms of $u_0$ and $v_0$, namely there exist $x_1,y_1, x_2, y_2 \in \mathbb{Z}$ such that:
$$u = x_1 u_0 + y_1 v_0 \quad \mbox{and} \quad v = x_2 u_0 + y_2 v_0.  $$

\noindent Then, $\Lambda' (u',v') = \{xu' + yv' =x(u+kv) + yv \ | \ x,y \in \mathbb{Z}\}$ with:
$$u' = (x_1 + k x_2) u_0 + (y_1 + k y_2) v_0  \quad \mbox{and} \quad  v' = v = x_2 u_0 + y_2 v_0. $$

\noindent Hence $\Lambda'$ is a periodic lattice for $d_\infty$ if and only if $u'$ is a translation vector of $d_\infty$, that is, $\Lambda' \subseteq \Lambda_0$, therefore, if and only if $u' \in \Lambda_0$, if and only if:
$$(x_1 + k x_2) \in \mathbb{Z}  \quad \mbox{and} \quad  (y_1 + k y_2) \in \mathbb{Z}.$$ 
The above conditions hold if and only if $k \in \mathbb{Z}$ (previous case of integral shearing) or $k \in \frac{1}{gcd(x_2,y_2)} \mathbb{Z}$. These are the only rational numbers satisfying the above conditions, and these conditions may change the lattice Bravais type. Any other values of $k$ are impossible since the only case that would satisfy the above conditions is $x_2 = y_2= 0$ and thus $v=0$, which contradicts the fact that $v$ is a periodic vector of $d_\infty$. 

\smallbreak
The analysis is analogous for a shearing of type $L_k$.

\begin{remark} \label{rem:shearing-identification}
It is interesting to consider the effect of a lattice $\mathbb{Z}$-shearing on the motif level. Using the same notation as above, the identification space of the flat motif diagram $d'$ gives rise to a motif diagram in the torus which is related to $d$ by a Dehn twist. A Dehn twist of a motif diagram $d$ induces a shearing of the corresponding flat motif diagram for the same longitude-meridian pair $(l,m)$ as $d$, and thus results in a sheared DP diagram. The same reasoning applies also to an admissible non-$\mathbb{Z}$ shearing of the basis, which gives rise to a new point lattice and thus to a sheared flat motif diagram and a sheared DP diagram. However, in this case of shearing  the two motif diagrams are not related by a Dehn twist. In fact this shearing is not visible on the motif level because of identification restrictions of the flat torus. These points are further addressed in Case 2 below. 
\end{remark}

\bigbreak

\noindent \textit{Case 2: Shearing of the DP diagram keeping the basis fixed.} 
\smallbreak

\noindent In this case we shear the DP diagram $d_\infty$ using a shearing transformation of the plane, but we retain the basis $B=(u,v)$ of the lattice of $d_\infty$ fixed. We shall first consider $\mathbb{Z}$-shearings of  $d_\infty$ by a matrix of type $M_k$ or $L_k$, with $k\in \mathbb{Z}$, as defined in~(\ref{shearing-matrices}). Keeping the basis $B=(u,v)$ fixed, and thus the lattice $\Lambda(u,v)$ point-wise fixed, this operation gives rise to a new DP diagram $(d'_\infty,\Lambda)$ and thus a new DP tangle $(\tau'_\infty,\Lambda)$. This was originally pointed out in \cites{Grishanov1,Grishanov.part1} and \cite{Farb}.  

\begin{figure}[H]
\centerline{\includegraphics[width=4in]{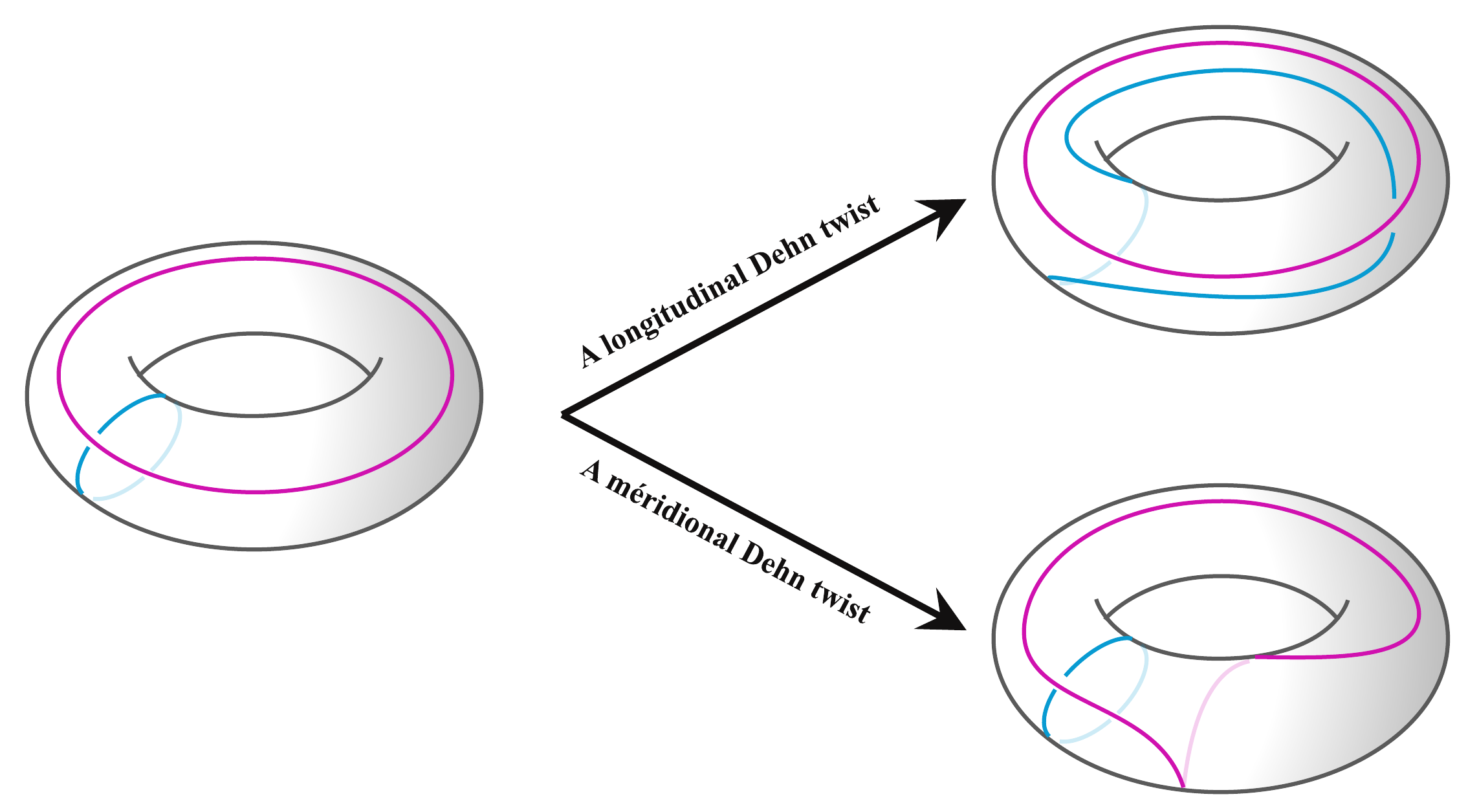}}
\vspace*{8pt}
\caption{\label{Rtwists} 
A longitudinal and a meridional Dehn twist of the torus.}
\end{figure}

Indeed, on the level of motif diagrams, a $\mathbb{Z}$-shearing of the DP diagram with a fixed basis corresponds to a Dehn twist of the torus, keeping the  longitude-meridian pair fixed. A Dehn twist gives rise to an orientation preserving self-homeomorphism of the torus, which is no more an isotopy, and it represents a full twist around a meridian or longitude. By convention, the identity matrix corresponds to the trivial Dehn twist. Figure~\ref{Rtwists} illustrates the effects of two types of Dehn twists corresponding to horizontal and  vertical  shears related to the matrices $L_k$ and $M_k$ respectively, with $k\in \mathbb{Z}$. 

\begin{figure}[th]
\centerline{\includegraphics[width=4in]{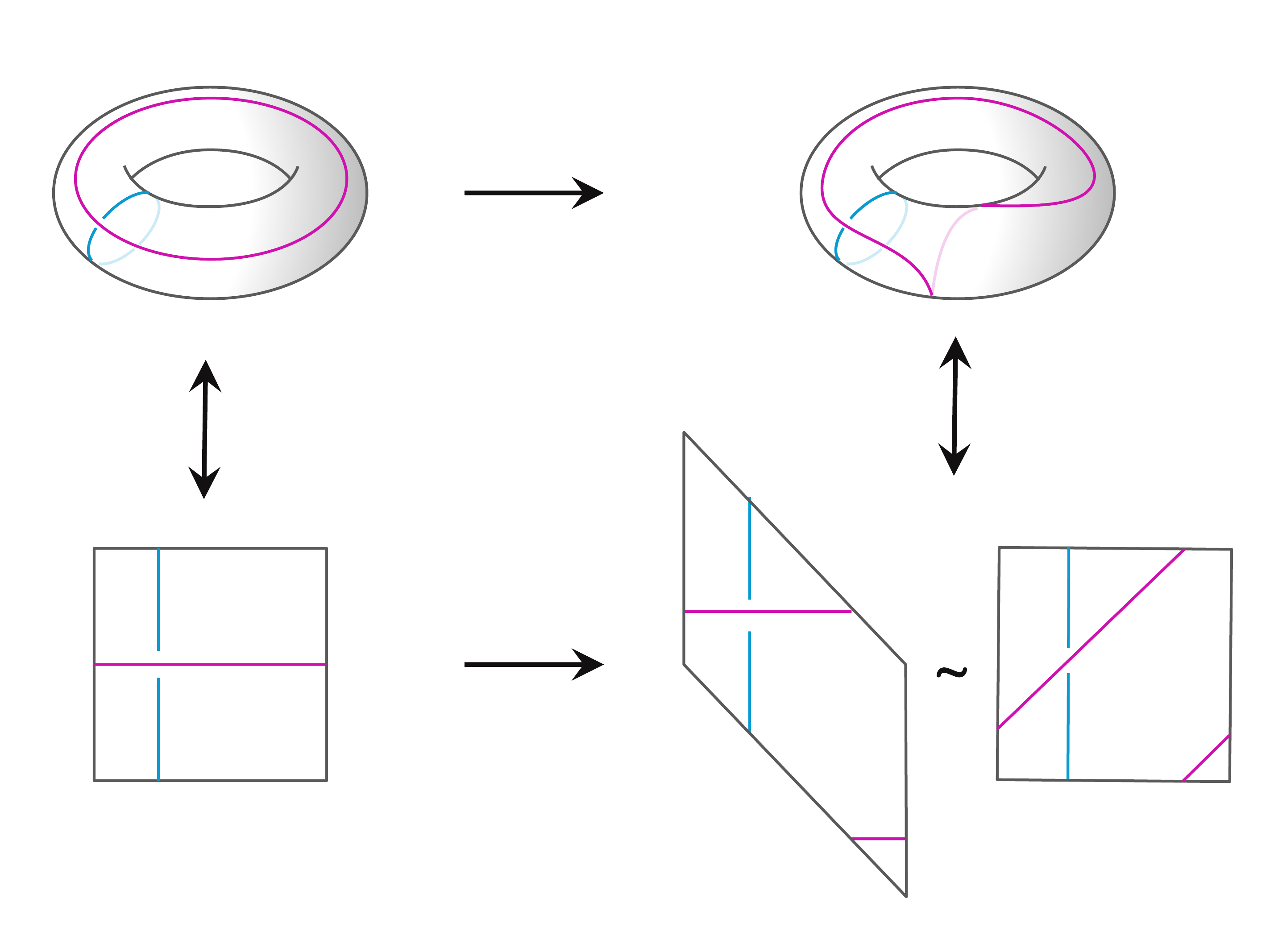}}
\vspace*{8pt}
\caption{\label{twist} 
A meridional Dehn twist on a flat motif, simultaneous or non-simultaneous.}
\end{figure}

Given that the torus $T^2= \mathbb{E}^2 / \Lambda$ arises as the identification space with respect to the longitude-meridian pair, this operation generates two different flat motif diagrams $d$ and $d'$, corresponding to the motif diagrams $(d,(l,m))$ and $(d',(l,m))$. Figure~\ref{twist} illustrates two motif diagrams related by a meridional Dehn twist and the induced transformation on the associated flat motif diagrams, keeping the basis and the lattice fixed (square flat torus).

In the universal cover, $d'$ gives rise to a new DP diagram $d^\prime_{\infty}$, which is a shearing of the original DP diagram $d_{\infty}$ for the same basis. The two corresponding sheared DP tangles are thus DP isotopic. The above procedure is exemplified in the right half of Figure~\ref{shearing}.

\bigbreak

Keeping the DP diagram fixed while shearing the basis by a shear factor $k \in \mathbb{Z}$, as addressed in Remark~\ref{rem:shearing-identification}, can be viewed as the reverse operation of shearing the DP diagram by a shear factor $-k \in \mathbb{Z}$  and keeping the basis fixed. This reverse operation can be generalized to non-integral shearings and all these connections are fully illustrated in Figure~\ref{shearing}. Intuitively, the sheared motif diagram $d'$ in the torus can be obtained by straightening (shearing) the parallelogram flat motif diagram $d$ into a rectangular one and identifying the opposite sides carrying along the diagram. 

With similar arguments as in Case 1, an arbitrary shearing of the DP diagram $(d_\infty,\Lambda)$ keeping the basis fixed and thus the lattice $\Lambda(u,v)$ point-wise fixed, will not, in general, give rise to a DP diagram $(d'_\infty,\Lambda)$, as the translational symmetry is in general not preserved. The only non-integral shearings preserving double periodicity are shearings by a matrix of type $M_k$ or $L_k$, with $k \in \frac{1}{gcd(x_2,y_2)} \mathbb{Z}$, as defined in~(\ref{shearing-matrices}). Such a  non-integral shearing of the DP diagram is justified as the reverse operation of an admissible non-integral shearing of the lattice by a factor $-k$ keeping the DP diagram fixed, in view also of Remark~\ref{rem:shearing-identification}.

\begin{definition}\label{def:dehn-equivalent}
Two (flat) motif diagrams of a DP tangle are said to be \textit{shear equivalent} if they are related by a Dehn twist of the torus or by an admissible non-$\mathbb{Z}$ shearing.
\end{definition}

\bigbreak

\noindent \textit{Case 3: Simultaneous shearing of the DP diagram and the basis.} 
\smallbreak

\noindent A simultaneous shearing of the DP diagram and its lattice $(d_\infty,\Lambda)$ by any matrix of type $M_k$ or $L_k$, as defined in~(\ref{shearing-matrices}), with arbitrary $k\in \mathbb{R}$,  results from a shearing of $\mathbb{E}^2$ applied to both $\Lambda(u,v)$ and $d_\infty$, that is, carrying along the DP diagram and its lattice. The resulting pair is a new DP diagram $(d'_\infty,\Lambda')$ and thus corresponds to a new DP tangle $(\tau'_\infty,\Lambda')$. This operation is visible on the level of the DP diagram, resulting in a corresponding shearing of the flat motif diagram. 
 On the level of motif diagrams, though, this operation is topologically trivial by the identification space of the torus, as illustrated on the left side of Figure~\ref{twist}. This generates in the torus the same motif diagram from different flat motif diagrams. 

\smallbreak
In conclusion, any shearing of a flat motif diagram, along with its flat torus, gives rise to a sheared DP diagram, but preserves the topology of the motif diagram.

\bigbreak

By the analysis of all three cases above and Definition~\ref{def:dehn-equivalent} we derive the following:

\begin{proposition}\label{prop:dehn-equivalent}
Two DP tangles are related by a shearing if and only if two (flat) motif diagrams of their are shear equivalent.
\end{proposition} 

Some remarks are now due.

\begin{remark} 
 A Dehn twist corresponds to an area-preserving transformation of the plane and thus it preserves minimality of a motif. So, applying any sequence of Dehn twists to a minimal motif will give rise to infinitely many shear equivalent minimal motifs.
\end{remark} 

\begin{remark}
With an eye on applying the above on materials modeled by DP tangles, it should be noted that there is a natural limitation on the amount of shearing a material can undergo due to its geometrical and physical properties. It would be very interesting to investigate this limitation for a given material, and this limitation would comprise a new geometrical invariant for the material. 
\end{remark}

\section{Unimodular linear transformations}\label{sec:unimodular}

A \textit{unimodular linear transformation of the plane} $\mathbb{E}^2$ is an area-preserving transformation that can be described by mapping one basis of the plane to another, determined by any matrix $M \in SL(2,\mathbb{R})$. The group $SL(2,\mathbb{R})$ is generated by the horizontal and vertical shearing matrices $M_k$ or $L_k$, as defined in~(\ref{shearing-matrices}) and analyzed in Section~\ref{sec:dehnequivalence}, cf. \cite[Lemma 8.1]{Lang1}. As another classical example, rotation matrices are also unimodular transformations, since $SO(2)$ is a subgroup of $SL(2,\mathbb{R})$. We recall that a unimodular transformation  $M=
\begin{pmatrix}
a & b \\
c & d 
\end{pmatrix}
\, (a, b, c, d \in \mathbb{R}, \ ad-bc=1)$ is classified by the absolute value of its trace, $|\mathrm{tr}(M)|=|a+d|$ (cf. \cite{Lang2} or \cite[Section~2.1]{Fushian}). If $|\mathrm{tr}(M)| <2$, $M$ is called \textit{elliptic}; if $|\mathrm{tr}(M)|=2$, $M$ is called \textit{parabolic}, and if $|\mathrm{tr}(M)|>2$,  $M$ is called \textit{hyperbolic}. 
It follows that the invariant curves for hyperbolic (resp. elliptic) linear transformations of $\mathbb{E}^2$ are hyperbolas (resp. ellipses), hence the terminology for these transformations. The parabolic transformations are so called by analogy, as intermediate between hyperbolic and elliptic. The fixed points are found by solving $z=\frac{az+b}{cz+d}$, and we see that a hyperbolic transformation has two fixed points in $\mathbb{R} \cup \{\infty\}$, one repulsive and one attractive, a parabolic transformation has one fixed point in $\mathbb{R} \cup \{\infty\}$, and an elliptic transformation has a pair of complex conjugate fixed points.

\bigbreak

\noindent \textit{Case 1: $|\mathrm{tr}(M)| < 2$.} 
\smallbreak

\noindent  If $|\mathrm{tr}(M)| < 2$, then $M \in SL(2,\mathbb{R})$ is said to be \textit{elliptic}.  A matrix in $SL(2,\mathbb{R})$ is  elliptic if and only if it is conjugate in $SL(2,\mathbb{R})$ to a unique rotation matrix $R \in SO(2)$: 
$R=
\begin{pmatrix}
\mathrm{cos} \, \theta & \mathrm{sin} \, \theta \\
-\mathrm{sin} \, \theta & \mathrm{cos}\, \theta 
\end{pmatrix}.
$ Note that  $\mathrm{tr}(R) = 2\, \mathrm{cos} \, \theta$.  
In $SL(2,\mathbb{Z})$, the only possibilities are matrices $M$ of finite order, namely matrices of orders $3$, $4$ and $6$, depending on the Bravais type of our lattice. On the DP tangle level, these transformations map $\tau_\infty$ either to a skewed version $\tau_\infty^3$, or to the same DP tangle $\tau_\infty^4$ rotated by $90^\circ$ clockwise or counterclockwise, or to another skewed version $\tau_\infty^6$, respectively. 
For the case of a matrix of order 4, $d$ is mapped to a motif $d'$, where the standard longitude and meridian of $d$ are mapped to the standard meridian and longitude of $d'$, respectively. For the case of a matrix $M$ of order 3 or 6, $d$ is mapped to a motif $d'$, where the standard longitude and meridian of $d$ are mapped to essential simple closed curves in the torus depending on $M$, respectively. 

\bigbreak

\noindent \textit{Case 2: $|\mathrm{tr}(M)| = 2$.} 
\smallbreak

\noindent  If $|\mathrm{tr}(M)| = 2$, then two cases must be distinguished in $SL(2,\mathbb{R})$ apart from the trivial case $M=I$. 
First, if $M = -I$, then $M$ represents a rotation of order 2, mapping $\tau_\infty$ to the same DP tangle $\tau_\infty^2$ rotated by $180^\circ$, clockwise or counterclockwise. 
On the torus level, the standard longitude and meridian of a motif $d$ are mapped to themselves in the image $d'$, up to orientation reversal. 
Second, every other matrix $M$ in this case is not diagonalizable, then $M$ is said to be \textit{parabolic} and maps $\tau_\infty$ to a sheared DP tangle $\tau_\infty'$, as detailed in the previous section. In this case fall the matrices of type $M_k$ or $L_k$, any product of matrices of type $M_k$, and any product of  matrices of type $L_k$, since all have trace equal to 2. A product $M_k \, L_r$ or $L_r  \, M_k$  falls in this case if $kr = 0$ (which is trivial) or $kr = -4$ for which there are various possibilities.
\bigbreak

\noindent \textit{Case 3: $|\mathrm{tr}(M)| > 2$.} 
\smallbreak

\noindent  If $|\mathrm{tr}(M)| > 2$, then $M$ is said to be \textit{hyperbolic}. A matrix in $SL(2,\mathbb{R})$ is hyperbolic if and only if it is diagonalizable over $\mathbb{R}$, or conjugate in $SL(2,\mathbb{R})$ to a unique matrix 
$
\begin{pmatrix}
a & 0 \\
0 & a^{-1} 
\end{pmatrix}
$
with $a \neq 1$. 
Any such matrix in $SL(2,\mathbb{R})$ maps $\tau_\infty$ to a stretched-compressed DP tangle $\tau_\infty'$, where one basic direction is expanded while the other is contracted, preserving the area. In $SL(2,\mathbb{Z})$, on the torus level, the standard longitude and meridian of a motif $d$ are mapped to two new simple closed curves of different slopes in the image $d'$, corresponding to the eigen-directions of $M$. These curves are preserved under iteration, but one is exponentially stretched while the other is exponentially contracted, yielding an Anosov action on the torus.

\section{DP tangle isotopy and flat motif equivalence} \label{sec:DPequivalent}

In this section we characterize DP tangle isotopy (as defined in  Definition~\ref{def:DPequivalence}) in the form of equivalence on the level of (flat) motif diagrams, called \textit{DP tangle equivalence}, which was the main target of the paper. We further make a discussion on minimal motifs.

\subsection{DP tangle equivalence}

In the previous sections, isotopy of DP tangles has been investigated thoroughly by highlighting two specific types of transformations: local isotopies and global isotopies. 

\smallbreak
Local isotopies of DP tangles are generated by planar isotopy moves and the three classical Reidemeister moves, preserving the lattice and the basis of $\mathbb{E}^2$, as detailed in Section~\ref{sec:local_isotopy}. On the level of (flat) motif diagrams, local isotopy is established using the theory of mixed links, thus obtaining Lemma~\ref{lem:localmotifisotopy}: \textit{Two (flat) motifs are isotopic if and only if two (flat) motif diagrams of theirs differ by a finite sequence of surface isotopy moves and the classical Reidemeister moves, as exemplified in Figures~\ref{Rmoves},~\ref{planar1},~\ref{planar2}, and ~\ref{planar3}.}

\smallbreak
Global isotopies of DP tangles are orientation-preserving affine transformations of the plane, without necessarily preserving the lattice and the basis. These are generated by translations (Section~\ref{sec:shiftequivalence}), scalings (Section~\ref{sec:scaleequivalence}) and shearings (Section~\ref{sec:dehnequivalence}), which generate all unimodular linear transformations (Section~\ref{sec:unimodular}). 

By a detailed analysis for each instance of global isotopies we translated the isotopy of DP tangles into an equivalence relation between their (flat) motif diagrams, after providing the appropriate definitions for (flat) motif diagrams (Definitions~\ref{def:shift-equivalent},~\ref{def:scale-equivalent},~\ref{def:dehn-equivalent}). For each transformation we distinguished  three cases: isotopy of the DP tangle keeping the lattice/basis fixed, transformation of the lattice/basis keeping the DP tangle fixed, and simultaneous transformation of both. We proved that the first two cases are equivalent to each other, while the third is invisible on the motif level. 

It follows from Propositions~\ref{prop:localmotifisotopy} ,~\ref{prop:shift-equivalent},~\ref{prop:scale-equivalent},~\ref{prop:dehn-equivalent} that isotopy of DP tangles can be characterized by any finite sequence of these local and global transformations on the motif level. The above analysis of the different instances of isotopy between two DP tangles lead to the following main result of our paper.

\begin{theorem}\label{th:equivalence}
Two DP tangles are DP isotopic if and only if two (flat) motif diagrams of theirs are related by a finite sequence of surface isotopy moves, Reidemeister moves, shift equivalence, scale equivalence, torus inflation/deflation, shear equivalence.
\end{theorem}

\begin{figure}[H]
\centerline{\includegraphics[width=5.6in]{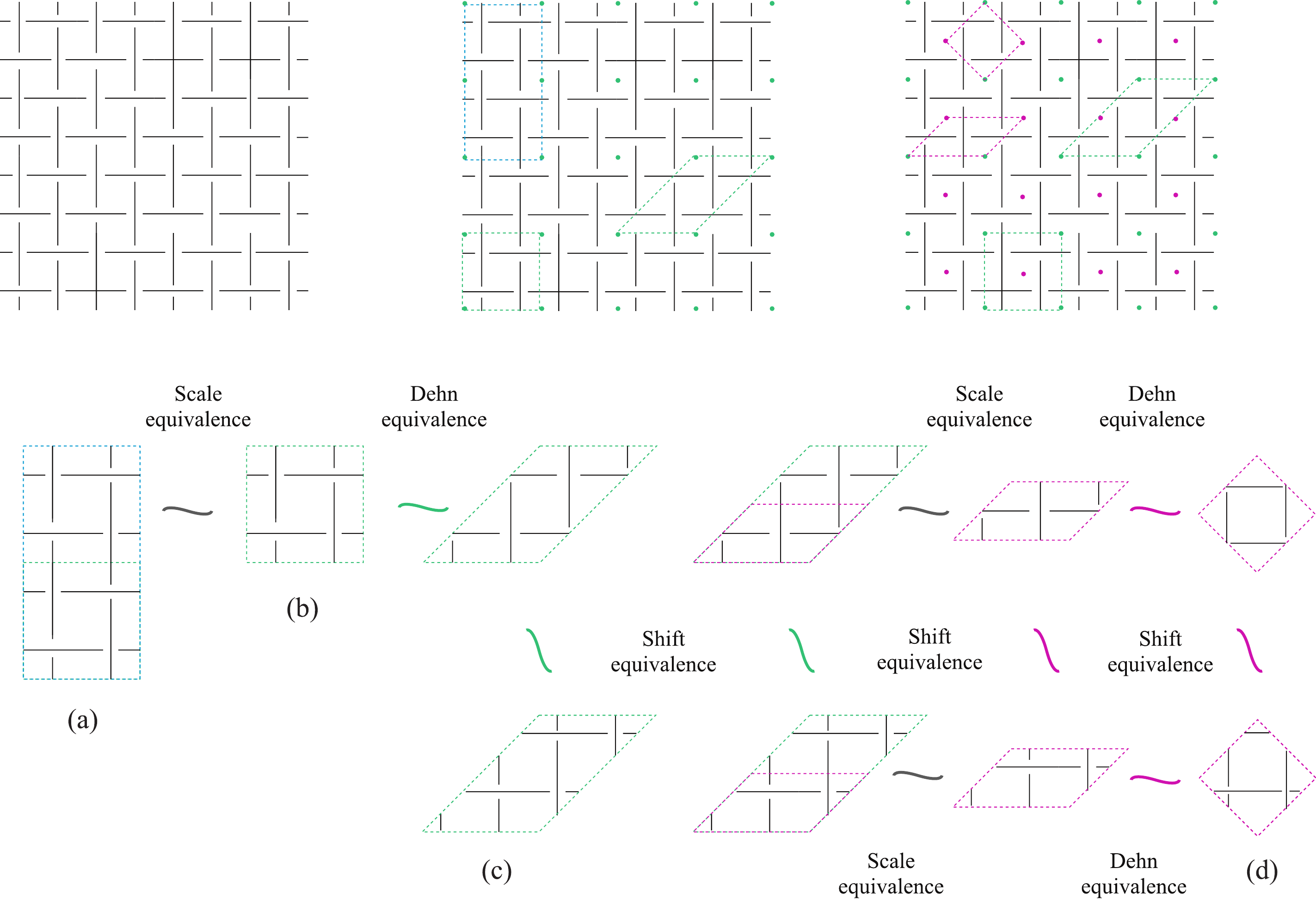}}
\vspace*{8pt}
\caption{\label{crossing} Different generating flat motifs for the same DP diagram and transformations among them culminating to the minimal motif (d).}
\end{figure}

Figure~\ref{crossing} captures many instances of Theorem~\ref{th:equivalence}, namely: shift equivalence, shear (Dehn) equivalence and scale equivalence. 

\subsection{A discussion on minimal motifs} \label{sec:minimalmotif}

We recall from Subsection~\ref{sec:lattice} the existence of a maximal periodic lattice for a DP tangle $\tau_\infty$, which corresponds to a minimal (flat) motif $\tau_{min}$ (Definition~\ref{def:minimal lattice}).

Finding a maximal lattice or, correspondingly, a minimal motif for a DP tangle can be very complex, as also noted  \cite{Grishanov1}. In \cite{Sonia1}, a method to find a minimal motif for a specific class of DP tangles, called untwisted weaves, is introduced by the third author, based on the homotopy classes of curves in the torus and combinatorial arguments. However, to our knowledge, there is no  general algorithm for finding a minimal motif of an arbitrary DP tangle. 

\begin{figure}[H]
\centerline{\includegraphics[width=5.3in]{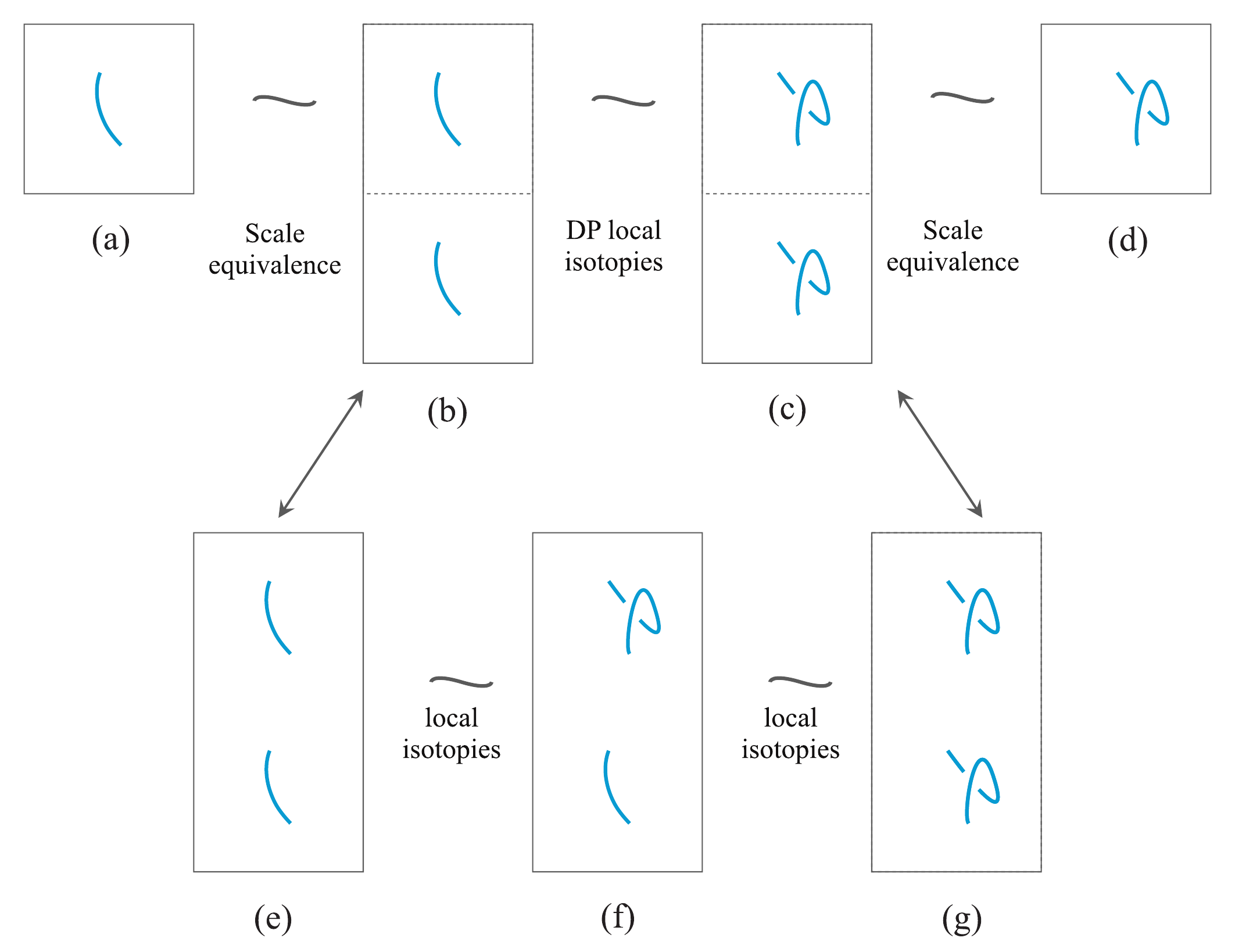}}
\caption{\label{R-move-minimal} 
A subtle obstruction to motif minimality.}
\end{figure}

\smallbreak
In Figure~\ref{R-move-minimal} we illustrate a subtle example. In the figure, we start with a minimal flat motif diagram (a) and we perform an R1 move, obtaining a second minimal flat motif diagram (d). Then we take a new flat motif diagram (e), which arises as the double cover (b) of (a), so this flat motif diagram is no longer minimal. On (e) we perform just one R1 move, obtaining the flat motif diagram (f), which is now minimal. If we performed both R1 moves in (e), we obtain (g) which arises as the double cover (c) of (d), and which in turn is no longer minimal. 

Another subtle situation is presented in Figure~\ref{crossing}, where from a flat motif diagram (b) (which could be mistaken to be minimal) we obtain an actual minimal flat motif diagram (d), via scale equivalence, Dehn twists and shift equivalence. Indeed, by Definition~\ref{def:minimal lattice} of a minimal motif, an arbitrary (flat) motif diagram $(d,(l,m))$ is related to a minimal (flat) motif diagram $(d_{min},(l',m'))$ by a finite covering map which may involve shift equivalence, Dehn twists and scale equivalence.

\section{DP tangle equivalence in other diagrammatic settings} \label{sec:othersettings}

In this section we discuss DP tangle equivalence when the links in the thickened torus, and subsequently their isotopy equivalences, belong to other diagrammatic categories. We first discuss the case of regular isotopy and framed link isotopy. Then we discuss the cases of virtual and welded links, the cases of pseudo and singular links, and the cases of tied and bonded links, all forming new classes of tangled objects. The theories of DP tangles in these other settings can find  applications in areas of science such as, molecular chemistry, textiles and meta-materials. 

\smallbreak
The general setting for DP tangles described in Section~\ref{sec:setup} and the detailed analysis in the subsequent sections, in particular the analysis of global transformations in Sections~\ref{sec:shiftequivalence},~\ref{sec:scaleequivalence},~\ref{sec:dehnequivalence},~\ref{sec:unimodular}, apply equally to any diagrammatic category. So, for the study of DP tangles related to any one of the above topological settings, we only need to adapt our analysis in the context of  local isotopy,  detailed in Section~\ref{sec:local_isotopy}. We present below each setting separately, calling especially on its combinatorial isotopy. The main  focus here lies in the elements of the diagrams which cross the specified longitude-meridian curves, so as to extend the surface isotopies of Section~\ref{sec:local_isotopy} to each diagrammatic category. Then, using the complete set of surface isotopies, any diagrammatic isotopy move in the category that crosses the longitude-meridian curves, can be pushed in the interior of a flat motif.

\subsection{DP tangle equivalence for regular isotopy and framed isotopy} \label{sec:DPregularframed}

\textit{Regular isotopy} is the equivalence relation between classical link diagrams generated by planar isotopy and only R2 and R3 Reidemeister moves. The move R1 is not allowed in this diagrammatic theory. Regular isotopy is introduced in \cite{Kauffman1987}, see also \cite{Kauffman1990}, where the Kauffman bracket polynomial is constructed as the regular isotopy equivalent of the Jones polynomial invariant for knots and links. An extension of the bracket polynomial for DP tangles has been constructed in \cites{Grishanov1,Grishanov.part2}. Regular isotopy projects to regular homotopy of the link projection, equivalently it preserves the Whitney degree of the link diagram. Finally, it has a natural interpretation when considering the link components as flat ribbons. 

\textit{Framed isotopy} is the topological equivalence of framed knots and links, see  \cite{Kauffman1990}. A framed knot is a knot endowed with a unit normal vector, hence it can be viewed as an embedded solid torus or, equivalently, an embedded annulus (a `ribbon') in 3-space, whereas a framed link is an embedding of one or more copies of a solid torus or an annulus. A framing unit (positive or negative) is, roughly speaking, a full twist of the solid torus or the ribbon. 

\begin{figure}[H]
\centerline{\includegraphics[width=2.5in]{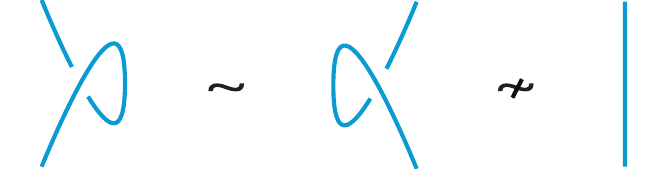}}
\caption{Framed R1 move.}
\label{Rframed}
\end{figure}

Representing a framed link by the central curves of its components, a framing unit is represented by a kink of the curve. So, classical links can represent framed links and the framing is registered as kinks in the diagram. There are two ways of projecting a framing unit  on the plane, as exemplified in Figure~\ref{Rframed}. So, this move is included in the diagrammatic \textit{framed isotopy}, along with planar isotopy and the Reidemeister moves R2 and R3. Clearly, R1 is excluded from  framed isotopy, as indicated in Figure~\ref{Rframed}. 

It is worth noting that regular isotopy on the surface of the sphere $S^2$ is equivalent to diagrammatic isotopy of framed links, as represented by classical links. This is not the case for the plane, where the framed R1 move cannot be realized by regular isotopy.

Framed links are used in the construction of closed, connected, orientable (c.c.o.) 3-manifolds via the surgery technique. The Witten invariant is a homeomorphism invariant of c.c.o. 3-manifolds based on the Jones polynomial. 
\bigbreak

Let now $\tau$ be a motif in the thickened torus and $d$ a diagram of $\tau$ undergoing regular isotopy. This isotopy in the torus is generated by the same local moves described above. This is periodically transmitted to the DP diagram $d_{\infty}$  of the DP tangle $\tau_{\infty}$. On the level of flat motif diagrams, the surface isotopy moves are valid (recall Figures~\ref{planar1},~\ref{planar2}, and~\ref{planar3}), so Proposition~\ref{prop:localmotifisotopy} adapts as follows:

\begin{proposition}\label{thm:regularmotifisotopy}
Two (flat) motif diagrams are \textit{regular isotopic} if and only if they differ by a finite sequence of surface isotopy moves and the classical Reidemeister moves R2 and R3.
\end{proposition}

\begin{figure}[ht]
\centerline{\includegraphics[width=5.7in]{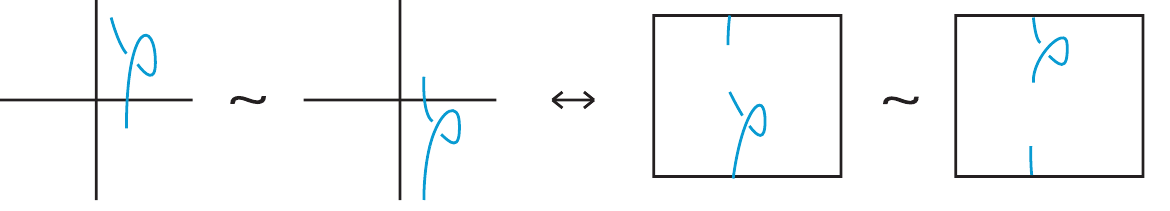}}
\caption{Framed surface isotopy move: a framing unit crossing the motif boundary.}
\label{planarframed}
\end{figure}

Suppose now that the diagram $d$ undergoes framed isotopy, then so does the DP diagram $d_{\infty}$. On the level of flat motif diagrams, we want to extend the surface isotopies by the \textit{framed surface isotopy move} indicating the passing of a framing unit through the meridian or longitude, as abstracted in Figure~\ref{planarframed}, where a framing unit is represented by a kink. Proposition~\ref{prop:localmotifisotopy} then adapts as follows:

\begin{proposition}\label{thm:framedmotifisotopy}
Two (flat) motif diagrams are \textit{framed isotopic} if and only if they differ by a finite sequence of  surface isotopy and framed surface isotopy moves, the framed R1 move and the Reidemeister moves R2 and R3.
\end{proposition}

The moves in Proposition~\ref{thm:framedmotifisotopy} generate the \textit{framed motif isotopy}. Furthermore, by the discussion in the beginning of this section, Theorem~\ref{th:equivalence}  adapts to the regular/framed isotopy setting, as follows:

\begin{theorem}\label{th:regularframedequivalence}
Two DP diagrams are \textit{DP regular isotopic} (resp. \textit{DP framed isotopic}) if and only if two (flat) motif diagrams of theirs are related by a finite sequence of regular isotopy moves (resp. framed motif isotopy moves), shift equivalence, scale equivalence, torus inflation/deflation, shear equivalence.
\end{theorem}

\subsection{Virtual and welded DP tangle equivalence } \label{sec:DPvirtualwelded} 

\begin{figure}[H]
\centerline{\includegraphics[width=4in]{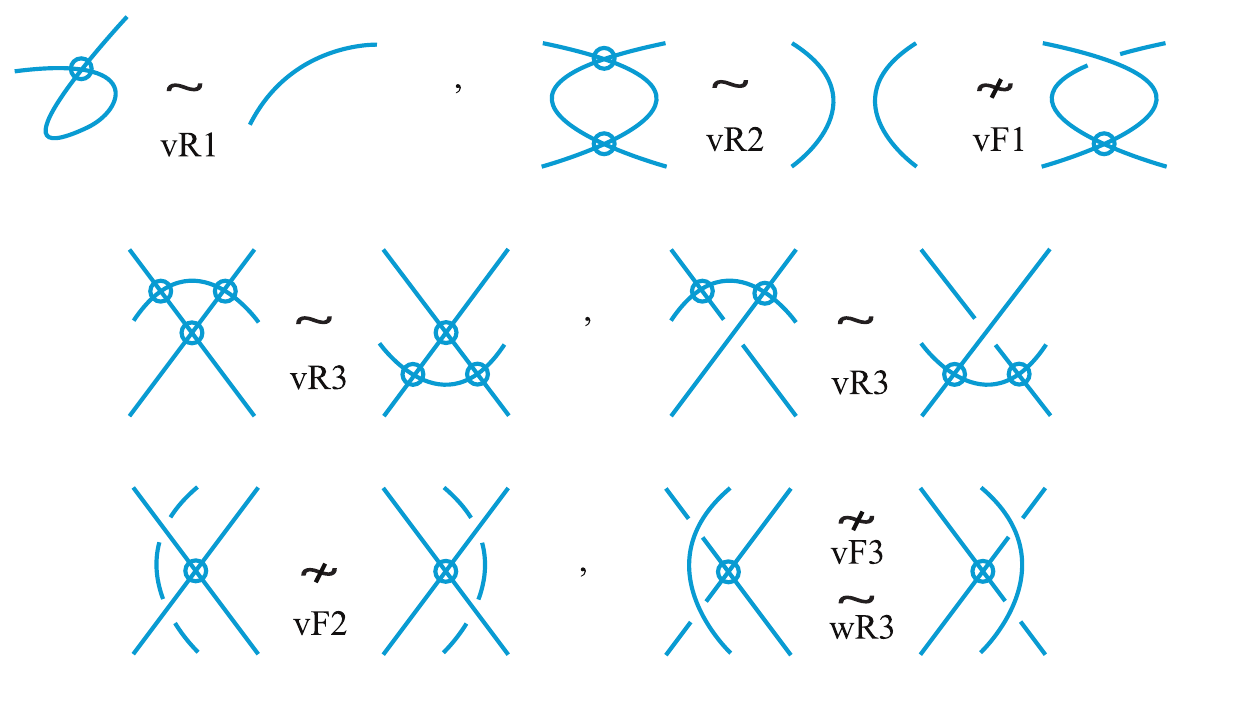}}
\caption{Virtual and welded moves: allowed and forbidden.}
\label{vmoves}
\end{figure}

Virtual knot theory is a diagrammatic extension of classical  knot theory in the sense that some crossings in a link diagram may represent just the permutation of the two arcs involved, with no further information of `over' or `under', and they are called \textit{virtual crossings}. The theory was introduced by Louis~H. Kauffman \cite{Kau2}. The diagrammatic equivalence in virtual knot theory, called \textit{virtual isotopy}, includes planar isotopy and the classical Reidemeister moves, extended by moves that contain virtual crossings. These moves are exemplified by the moves vR1, vR2, vR3 in Figure~\ref{vmoves} (where also the variant with the opposite type of classical crossing is assumed) and they are all special cases of the universal `detour move' whereby an arc containing all virtual crossings can slide across any parts of the diagram. In this theory we also have the virtual forbidden moves vF1, vF2 and vF3 depicted in Figure~\ref{vmoves}.

\begin{figure}[H]
\centerline{\includegraphics[width=4in]{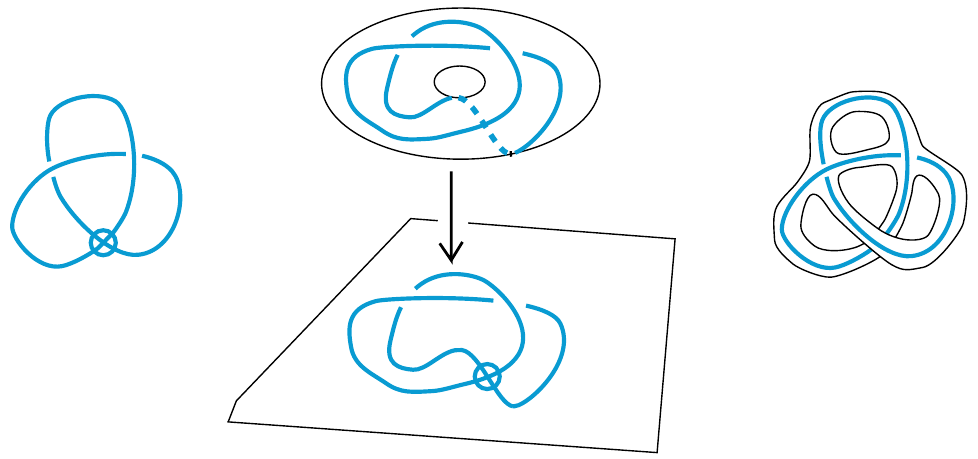}}
\caption{Surface realizations of virtual knots.}
\label{virtual1}
\end{figure}

Virtual links have an interesting interpretation as embeddings of  links in thickened surfaces, taken up to addition and subtraction of empty handles \cite{CKS}:  a  virtual crossing is regarded as  a detour of one of the arcs in the crossing through a 1-handle that has been attached to the 2-sphere of the original diagram  (see middle illustration of Figure~\ref{virtual1}). Another nice interpretation of a virtual link diagram is obtained by forming a ribbon--neighborhood surface of the diagram, where a virtual crossing is represented by abstract ribbons passing over one another without interacting \cite{KamadaNS}, as in the right part of Figure~\ref{virtual1}. This consideration could find a physical application in materials science, where the prevention of friction between strands of an embedded DP tangle may be interesting to consider when defining the energy/relaxing state of the corresponding material.  

A  refinement of virtual knot theory is \textit{welded knot theory}, introduced in \cite{Fenn}. In this diagrammatic theory the moves generating virtual isotopy are all included but there is additionally the wR3 move (which is the forbidden move vR3 in virtual knot theory), as depicted in  Figure~\ref{vmoves}. This move contains an over arc and one virtual crossing; in general it enables to detour sequences of classical crossings \textit{over} welded crossings. Hence, welded knot theory can be realized as the quotient of virtual knot theory  modulo the  wR3  move. The explanation for the choice of moves lies in the fact that the move  wR3 preserves the combinatorial fundamental group. This is not the case for the other forbidden move vF2, so it remains forbidden also for welded links. 

\begin{figure}[H]
\centerline{\includegraphics[width=5in]{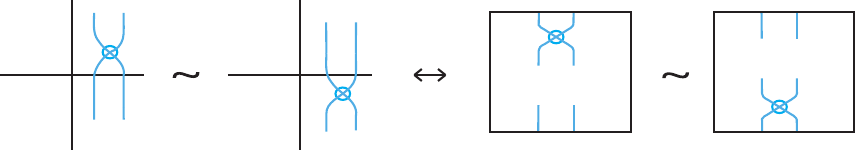}}
\caption{Virtual and welded motif isotopy moves of type e).}
\label{vsurfacemoves}
\end{figure}

Let now $d$ be a virtual (resp. welded) motif diagram in the torus undergoing virtual (resp. welded) isotopy. This isotopy in the torus is generated by the same local moves described above. This is periodically transmitted to the DP diagram $d_{\infty}$. On the level of flat motifs, the \textit{virtual (resp. welded) surface isotopy move}, whereby a virtual (resp. welded) crossing passes through the meridian or longitude (see Figure~\ref{vsurfacemoves}), extends the usual surface isotopy moves (recall Figures~\ref{planar1},~\ref{planar2}, and~\ref{planar3}).

So, as stated above, any one of the virtual isotopy moves on a motif diagram, which crosses the specified meridian or longitude on the torus, can be pushed to the interior of the flat motif diagram. Hence, Proposition~\ref{prop:localmotifisotopy} adapts as follows: 

\begin{proposition}\label{thm:virtual/weldedmotifisotopy}
Two virtual (resp. welded) (flat) motif diagrams are \textit{virtual (resp. welded) isotopic} if and only if they differ by a finite sequence of surface isotopy moves, virtual (resp. welded) surface isotopy moves, the classical Reidemeister moves and the (allowed) virtual (resp. welded) moves.
\end{proposition} 

The moves in Proposition~\ref{thm:virtual/weldedmotifisotopy} generate the \textit{virtual (resp. welded) motif isotopy}. 

\begin{remark}
In the context of DP tangles, a first observation is that any motif, being by definition a link in the thickened torus, can be represented by a virtual link diagram, as in the bottom and left part of Figure~\ref{virtual1}.  In \cite{Kawauchi}, particular cases of DP diagrams have been studied, which do not contain closed components. Equivalence of the corresponding unit-square flat motif diagrams, called \textit{$(m,n)$-knitting patterns}, is considered up to local isotopy preserving the boundary of the square. Then, by considering the virtual link diagram associated to any such flat motif diagram, a characterization of equivalence of flat motif diagrams is stated up to equivalence of the corresponding virtual link diagrams. 
\end{remark}

Moreover, by the discussion in the beginning of this section, Theorem~\ref{th:equivalence} adapts for virtual (resp. welded) DP diagrams, as follows:

\begin{theorem}\label{th:virtualweldedequivalence}
Two virtual (resp. welded) DP diagrams are \textit{DP virtual isotopic} (resp. \textit{DP welded isotopic}) if and only if two (flat) motif diagrams of theirs are related by a finite sequence of virtual (resp. welded) motif isotopy moves, shift equivalence, scale equivalence, torus inflation/deflation, shear equivalence.
\end{theorem}

\subsection{Singular and pseudo DP tangle equivalence } \label{sec:DPpseudosingular} 

\begin{figure}[ht]
\centerline{\includegraphics[width=3.5in]{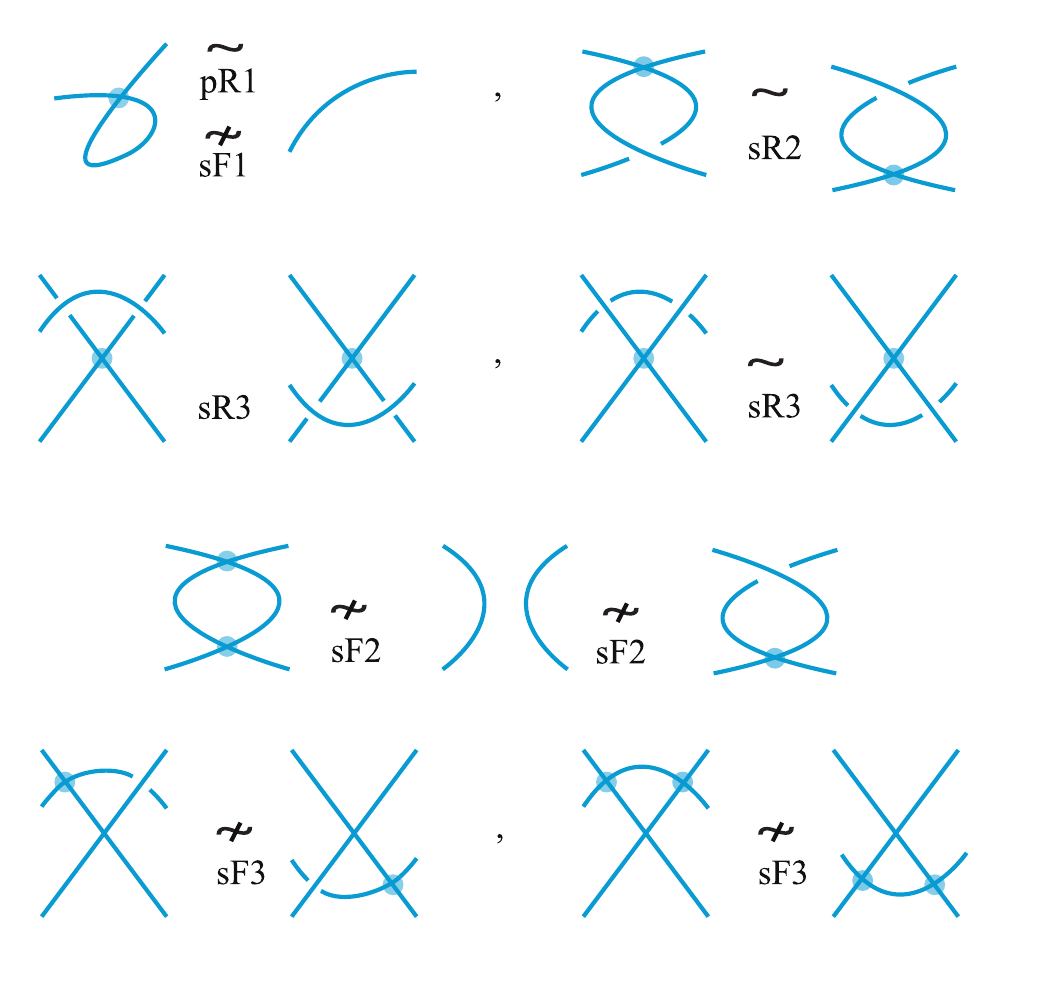}}
\vspace*{8pt}
\caption{Singular and pseudo-knot moves: allowed and forbidden.}
\label{Rsingularpseudo}
\end{figure}

\textit{Singular knot theory} appeared in the context of the theory of Vassiliev's finite type knot invariants (see references in \cites{Grishanov.Vassiliev1,Grishanov.Vassiliev2}). Singular knots are knots with finitely many rigid self-intersections, the singular crossings, interpreted as rigid vertices in a spatial graph. So singular link isotopy includes classical link isotopy together with rigid vertex isotopy. Figure~\ref{Rsingularpseudo} exemplifies the diagrammatic moves in the theory that extend planar isotopy and the classical Reidemeister moves, as well as the  singular forbidden moves, sF1, sF2 and sF3 (where the middle crossing could be classical or singular). 

Another related diagrammatic category is the  \textit{theory of pseudo knots}. Pseudo diagrams of knots, links and spatial graphs were introduce by Hanaki in \cite{H} as projections on the 2-sphere with apparently missing,  rather hidden,  crossing information on some crossings, called \textit{precrossings}. 

The theory of singular knots is closely related to the theory of pseudo knots. Namely, the diagrammatic \textit{pseudo link isotopy} is generated by planar isotopy and the classical Reidemeister moves extended by the singular isotopy moves and the pseudo-Reidemeister move 1, denoted pR1, which is the forbidden move sF1 in singular knot theory (all exemplified in Figure~\ref{Rsingularpseudo}). Hence, pseudo knot theory can be realized as the quotient of singular knot theory, modulo the pR1 move. 

Pseudo knots make up a novel and significant model for DNA knots since there are some DNA knots where it is difficult to distinguish (even with electron microscope) the relative positions of the two arcs in some crossings. Analogous situation can certainly occur also in \textit{worn textiles}, which thus can be modeled by pseudo DP tangles.    

\begin{figure}[H]
\centerline{\includegraphics[width=5in]{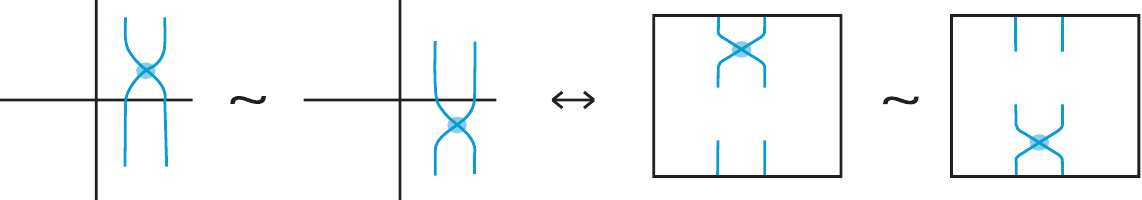}}
\caption{Singular and pseudo motif isotopy moves of type e).}
\label{spsurfacemoves}
\end{figure}

\bigbreak
Let now $\tau$ be a singular motif in the thickened torus and $d$ a diagram of $\tau$ undergoing singular isotopy. Respectively, let  $d$ be a motif diagram in the torus undergoing  pseudo knot isotopy. These isotopies in the torus are generated by the same local moves described above. Each isotopy carries through periodically to the DP diagram $d_{\infty}$  of the singular resp. pseudo DP tangle $\tau_{\infty}$. On the level of flat motif diagrams, the \textit{singular (resp. pseudo) surface isotopy move}, whereby a singular crossing (resp. a precrossing) passes through the meridian or longitude (see Figure~\ref{spsurfacemoves}), extends the usual surface isotopy moves (recall Figures~\ref{planar1},~\ref{planar2}, and~\ref{planar3}).

So, Proposition~\ref{prop:localmotifisotopy} adapts as follows: 

\begin{proposition}\label{thm:singular/pseudomotifisotopy}
Two singular (resp. pseudo) (flat) motif diagrams are \textit{singular (resp. pseudo) isotopic} if and only if they differ by a finite sequence of surface isotopy moves, singular (resp. pseudo) surface isotopy moves, the classical Reidemeister moves and the (allowed) singular (resp. pseudo) isotopy moves.
\end{proposition} 

The moves in Proposition~\ref{thm:singular/pseudomotifisotopy} generate the \textit{singular (resp. pseudo) motif isotopy}. Moreover, by the discussion in the beginning of this section, Theorem~\ref{th:equivalence}  adapts for singular  DP tangles, as follows:

\begin{theorem}\label{th:singularequivalence}
Two singular DP diagrams are \textit{DP singular isotopic} if and only if two (flat) motif diagrams of theirs are related by a finite sequence of singular motif isotopy moves, shift equivalence, scale equivalence, torus inflation/deflation, and shear equivalence.
\end{theorem}

Using the theory of singular motifs, finite type invariants have been constructed in \cites{Grishanov.Vassiliev1,Grishanov.Vassiliev2} to distinguish DP tangles which were not differentiated by the bracket-type polynomial invariant constructed in \cites{Grishanov1,Grishanov.part2}. The theory of singular DP tangles could also find interesting applications in molecular chemistry. 

\bigbreak

The theory of pseudo DP tangles and their equivalence has been treated in \cite{DLM-pseudo}. By  the discussion in the beginning of this section, Theorem~\ref{th:equivalence}  adapts for pseudo  DP tangles as well. However, here there is a subtlety related to the scale equivalence, due to the fact that for a precrossing we do not know whether it an over or an under crossing, yet the information is there. Subsequently, in \cite{DLM-pseudo} we needed to introduce the notion of \textit{homologous pre-crossings}, which are precrossings bearing the same (hidden) information. Hence, two motifs can be scale equivalent only if the corresponding pre-crossings that they contain are homologous. So, Theorem~\ref{th:equivalence} adapts as follows:

\begin{theorem}\label{th:pseudoequivalence}
Two singular (resp. pseudo) DP diagrams are \textit{DP singular isotopic} (resp. \textit{DP pseudo isotopic}) if and only if two (flat) motif diagrams of theirs are related by a finite sequence of singular (resp. pseudo) motif isotopy moves, shift equivalence, scale equivalence, torus inflation/deflation, shear equivalence.
\end{theorem}

\subsection{Tied and bonded DP tangle equivalence } \label{sec:DPtied} 

\begin{figure}[H] 
\centerline{\includegraphics[width=3.5in]{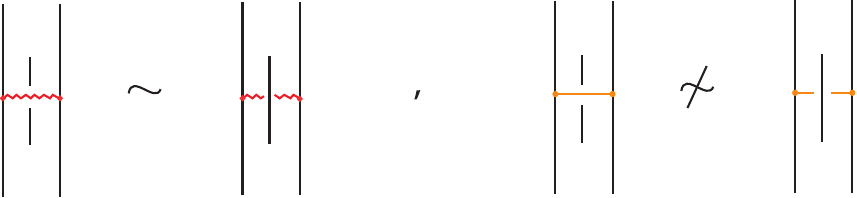}} 
\caption{The tied elementary move, which is forbidden for bonded links.} 
\label{tbmove} 
\end{figure}

\begin{figure}[H]
\begin{center}
\includegraphics[width=3.5in]{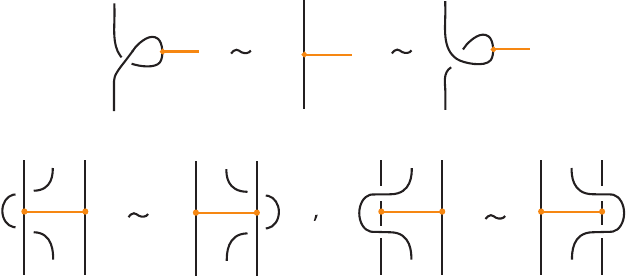}
\end{center}
\caption{Topological vertex isotopy moves for bonded and tied links.}
\label{tightiso_top}
\end{figure}

\textit{Tied links} were introduced in \cite{AJ1} as generalization of links in $S^3$. A tied link is  a classical link  equipped with `ties'. A tie connects two points of the link and behaves like a phantom: its ends can slide along the arcs that it connects, passing across any other parts of the link without obstruction, see left-hand illustration of Figure~\ref{tbmove}, where a tie is depicted as red spring. There are further the rules that: two ties joining the same pair of components merge into one tie, and ties joining points of the same component can be deleted or introduced at will. Hence, a set of ties in a link diagram provides a combinatorial structure for defining a partition of the set of components of the link, by considering two components of the tied link to belong to the same partition subset if there is a tie connecting them. Then, \textit{tied link isotopy} is defined as ambient isotopy between the underlying links (ignoring the ties), taking also into consideration that the set of ties in the links define the same partition of the set of components. In terms of diagrams we have planar isotopy and the Reidemeister moves together with local moves involving ties, as illustrated in Figures~\ref{tbmove} and \ref{tightiso_top}, where the orange springs represent ties.

It is worth adding that the above results generalize directly for the diagrammatic theory of \textit{tied singular links}, introduced and studied in \cite{AJ3} as a generalization of singular links, as well as for \textit{tied pseudo links}, introduced and studied in \cite{D1}.  

\bigbreak
Considering  embedded simple arcs in $S^3$ in place of the ties, we obtain the theory of \textit{bonded links}. More precisely, a bonded link is a pair $(L, B)$, where $L$ is a link in $S^3$ and $B$ is a set of $k$ pairwise disjoint  simple arcs, the `bonds',  properly embedded in the complement $S^3\backslash L$ of the link, such that the boundaries of the  bonds intersect the link transversely in $2k$ distinct points. The intersection points are not considered as rigid vertices. Bonds were introduced in \cite{GGLDSK} in the context of bonded knotoids, under rigid vertex isotopy, for modeling open knotted protein chains. 

\begin{figure}[H]
\begin{center}
\includegraphics[width=3.3in]{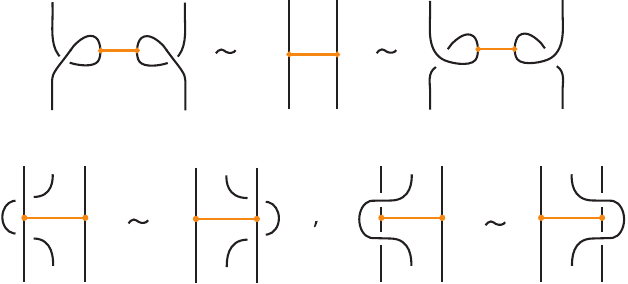}
\end{center}
\caption{Rigid  vertex isotopy moves for bonded links.}
\label{tightiso_rigid}
\end{figure}

Bonded link isotopy, either \textit{rigid vertex isotopy} or \textit{topological vertex isotopy} is defined as ambient isotopy between links, taking also into consideration the set of bonds in the links (cf. \cite{DLK} and references therein). In terms of diagrams we have planar isotopy and the Reidemeister moves together with local moves involving bonds, as illustrated in Figures~\ref{tightiso_top} and~\ref{tightiso_rigid}, where bonds are depicted as orange line segments. We also have forbidden moves in the theory, as for example the move depicted in the right-hand side of Figure~\ref{tbmove}. Hence, tied knot theory can be realized as the quotient of the topological vertex bonded knot theory, modulo the tied elementary move, Figure~\ref{tbmove} and the tie cancellations/additions.

\begin{figure}[H] 
\centerline{\includegraphics[width=5.7in]{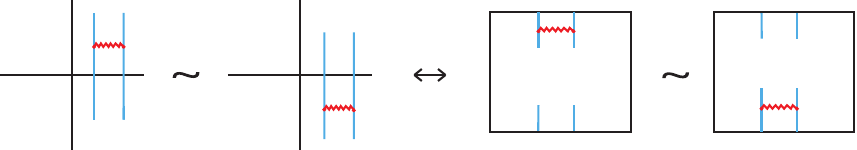}}
\caption{Bonded/tied surface isotopy moves.}
\label{tbmotifisotopy}
\end{figure}

Let now $\tau$ be a bonded  (resp. tied)   motif in the thickened torus and $d$ a diagram of $\tau$ undergoing bonded  (resp. tied)  isotopy. This isotopy is generated in the torus  by the same diagrammatic local moves described above.  Each isotopy carries through periodically to the DP diagram $d_{\infty}$  of the DP tangle $\tau_{\infty}$. On the level of flat motifs, the \textit{bonded  surface isotopy move} (resp. the \textit{tied surface isotopy move}), whereby a bond (resp. tie)   passes through the meridian or longitude (see Figure~\ref{tbmotifisotopy}), extends the usual surface isotopy moves  (recall Figures~\ref{planar1},~\ref{planar2}, and~\ref{planar3}).

So, Proposition~\ref{prop:localmotifisotopy} adapts as follows: 

\begin{proposition}\label{thm:tied/bondedmotifisotopy}
Two bonded (resp. tied) (flat) motif diagrams are \textit{bonded rigid vertex or topological vertex isotopic  (resp. tied isotopic)} if and only if they differ by a finite sequence of surface isotopy moves, bonded  (resp. tied)  surface isotopy moves, the classical Reidemeister moves and the respective bonded (resp. tied) isotopy moves.
\end{proposition} 

The moves in Proposition~\ref{thm:tied/bondedmotifisotopy} generate the \textit{bonded rigid vertex or topological vertex (resp. tied)  motif isotopy}. Moreover, by the discussion in the beginning of this section, Theorem~\ref{th:equivalence} applies for bonded  (resp. tied)  DP tangles, as follows:

\begin{theorem}\label{th:tiedbondedequivalence}
Two bonded (resp. tied) DP diagrams are \textit{DP bonded rigid vertex or topological vertex  isotopic} (resp. \textit{DP tied isotopic}) if and only if two (flat) motif diagrams of theirs are related by a finite sequence of bonded rigid vertex or topological vertex  (resp. tied) motif isotopy moves, shift equivalence, scale equivalence, torus inflation/deflation, shear equivalence.
\end{theorem}

\end{document}